\documentclass[11pt]{amsart}
\usepackage[hmargin=3cm,vmargin=3cm]{geometry}
\usepackage[backend=bibtex,doi=false,isbn=false,url=false,style=alphabetic,maxnames=6]{biblatex}
\usepackage{etex}
\bibliography{../../../bibliography/bib}

\usepackage[latin1]{inputenc}
\usepackage{amssymb,amsfonts}
\usepackage{amsthm,amsmath,amscd,stmaryrd,latexsym,youngtab,mathtools,pinlabel}
\usepackage[scaled]{helvet} 
\usepackage{courier} 

\usepackage{euscript,mathrsfs}

\usepackage[nohug]{diagrams}
\usepackage{tikz, tikz-cd}
\usetikzlibrary{matrix, arrows, calc, fit}

\RequirePackage{color}
\definecolor{myred}{rgb}{0.75,0,0}
\definecolor{mygreen}{rgb}{0,0.5,0}
\definecolor{myblue}{rgb}{0,0.25,0.65}
\definecolor{references}{rgb}{0,0,1}

  \RequirePackage[pdftex,
   colorlinks = true,
   urlcolor = references, 
   citecolor = references, 
   linkcolor = references, 
   ]
 {hyperref}

\newtheorem{theorem}{Theorem}[section]
\newtheorem{lemma}[theorem]{Lemma}

\newtheorem{proposition}[theorem]{Proposition}

\newtheorem{corollary}[theorem]{Corollary}

\newtheorem{slogan}[theorem]{Slogan}
\newtheorem{observation}[theorem]{Observation}

\theoremstyle{definition}

\newtheorem{definition}[theorem]{Definition}

\newtheorem{notation}[theorem]{Notation}

\newtheorem{example}[theorem]{Example}

\theoremstyle{remark}
\newtheorem{remark}[theorem]{Remark}

\newtheorem{question}[theorem]{Question}

\numberwithin{equation}{section}

\def\AB{{\mathbf A}}    \def\AC{{\mathcal{A}}}
\def\BB{{\mathbf B}}    \def\BC{{\mathcal{B}}}
\def\CB{{\mathbf C}}    \def\CC{{\mathcal{C}}}
\def\DB{{\mathbf D}}    
\def\EB{{\mathbf E}}    
\def\FB{{\mathbf F}}

\def\IB{{\mathbf I}}    \def\IC{{\mathcal{I}}}
    \def\JC{{\mathcal{J}}}
    \def\KC{{\mathcal{K}}}
    \def\LC{{\mathcal{L}}}
    \def\MC{{\mathcal{M}}}

\def\PB{{\mathbf P}}    
    
\def\RB{{\mathbf R}}

\def\TB{{\mathbf T}}    \def\TC{{\mathcal{T}}}
\def\UB{{\mathbf U}}    
\def\VB{{\mathbf V}}

\def\AS{{\EuScript A}}
\def\BS{{\EuScript B}}
\def\CS{{\EuScript C}}

\def\FS{{\EuScript F}}

\def\MS{{\EuScript M}}
\def\NS{{\EuScript N}}

\def\TS{{\EuScript T}}

\def\a{\alpha}
\def\b{\beta}
\def\g{\gamma}

\def\d{\delta}

\def\e{\varepsilon}
\def\k{\kappa}
\def\l{\lambda}

\let\phi=\varphi

\usepackage{bbm}

\def\N{{\mathbbm N}}
\def\R{{\mathbbm R}}
\def\Z{{\mathbbm Z}}

\def\1{\mathbbm{1}}

\renewcommand{\to}{\rightarrow}

\def\modules{{\textrm{-mod}}}
\newcommand{\proj}{\textrm{-proj}}

\def\bimodules{{\textbf{-bimod}}}

\newcommand{\refequal}[1]{\xy {\ar@{=}^{#1}
(-1,0)*{};(1,0)*{}};
\endxy}

\newcommand{\ip}[1]{\langle #1 \rangle}
\renewcommand{\matrix}[1]{\begin{bmatrix}#1\end{bmatrix}}
\newcommand{\inv}{^{-1}}
\renewcommand{\setminus}{\smallsetminus}
\renewcommand{\d}{\delta}
\renewcommand{\emptyset}{\varnothing}

\newcommand{\one}{\mathbbm{1}}

\newcommand{\smMatrix}[1]{\left[\begin{smallmatrix}#1\end{smallmatrix}\right]}

\newcommand{\Id}{\operatorname{Id}}
\newcommand{\Hom}{\operatorname{Hom}}

\newcommand{\im}{\operatorname{im}}
\newcommand{\Tot}{\operatorname{Tot}}

\newcommand{\Cone}{\operatorname{Cone}}

\newcommand{\End}{{\rm End}}

\newcommand{\Ext}{{\rm Ext}}
\newcommand{\Tor}{{\rm Tor}}

\newcommand{\op}{{\rm op}}

\newcommand{\Mat}{\operatorname{Mat}}
\newcommand{\Cob}{\operatorname{Cob}}
\renewcommand{\k}{\mathbbm{k}}
\newcommand{\idempotent}{\idemp}
\newcommand{\idemp}{\EB}
\newcommand{\otherIdemp}{\FB}
\newcommand{\otherIdempotent}{\otherIdemp}
\newcommand{\unital}{\UB}
\newcommand{\counital}{\CB}
\newcommand{\homIdempt}{\TB}
\newcommand{\Tate}{\homIdempt}
\newcommand{\cotimes}{\:\hat{\otimes}\:}

\newcommand{\K}{\mathcal{K}}

\newcommand{\Homb}{\Hom^\bullet}
\newcommand{\Homg}{\operatorname{hom}}
\newcommand{\Endg}{\operatorname{end}}

\newcommand{\hocolim}{\operatorname{hocolim}}

\begin{document}

\begin{abstract}
In these notes we develop some basic theory of idempotents in monoidal categories.   We introduce and study the notion of a pair of complementary idempotents in a triangulated monoidal category, as well as more general idempotent decompositions of identity.  If $\idemp$ is a categorical idempotent then $\End(\idemp)$ is a graded commutative algebra.  The same is true of $\Hom(\idemp,\idemp^c[1])$ under certain circumstances, where $\idemp^c$ is the complement.  These generalize the notions of cohomology and Tate cohomology of a finite dimensional Hopf algebra, respectively.
\end{abstract}

\title[Categorical idempotents]{Idempotents in triangulated monoidal categories }

\author{Matthew Hogancamp} \address{University of Southern California}\thanks{The author was partially supported by NSF grant DMS-1255334.}

\maketitle

\tableofcontents


\section{Introduction}\label{sec:intro}



%

Let $\AS$ be a monoidal category with monoid $\otimes$ and monoidal identity $\one$.  A \emph{unital idempotent} in $\AS$ is an object $\unital\in\mathcal{A}$ together with a morphism $\eta:\one \rightarrow \unital$ so that $\eta\otimes\Id_\unital: \one\otimes \unital\rightarrow \unital\otimes \unital$ and $\Id_\unital\otimes\eta: \unital\otimes \one \rightarrow \unital\otimes \unital$ are isomorphisms.  Similarly, a \emph{counital idempotent} in $\AS$ is an object $\counital\in \AS$ together with a morphism $\e:\counital\rightarrow \one$ such that $\e\otimes \Id_\counital:\counital\otimes \counital\rightarrow \one\otimes \counital$ and $\Id_\counital\otimes \e:\counital\otimes \counital\rightarrow \counital\otimes \one$ are isomorphisms.

Such objects are well-studied and appear naturally in a number of related contexts, for instance Bousfield localization \cite{Bous79_loc,KrauseNotes}, semi-orthogonal decompositions \cite{BonOrl-semi}, and recollement \cite{BBD81}.  See \S \ref{subsec:introDeomp}.

Our purpose in writing this note is to have a convenient reference for future work in categorification and higher representation theory, but we also hope that non-experts will find this to be an accessible (and elementary) introduction to this interesting and useful theory.  Other useful references are the preprint of Drinfeld-Boyarchenko \cite{BoyDrin-idemp} (a subset of which is is found in published form in \S 2 of \cite{BD14}) and also Balmer-Favi \cite{BalmerFavi11}.  The terminology of (co)unital idempotent is new, as far as we can tell; in \cite{BD14} they are called closed and open idempotents; in \cite{BalmerFavi11} they are called right and left idempotents; in topology literature (see for instance Section 4.2 of \cite{BorceuxBook2}) they are called idempotent monads and comonads, or sometimes localization and colocalization functors.  We prefer our terminology, because it is descriptive and context-independent.

\begin{remark}
The present document began its life as an appendix to \cite{Hog15}.  We decided to separate the two, after which point we added the material which now appears in \S \ref{sec:triangMonoidal}, \S \ref{sec:decompOfOne}, and \S \ref{sec:tate}.  
\end{remark}


\subsection{Familar examples}
Before continuing with the remainder of the introduction, we give some examples which the reader may wish to have in mind:
\begin{enumerate}
\item  There $\Z\rightarrow \Z_{(p)}$ which makes ``$\Z$ localized at $p$'' into a unital idempotent in the category of abelian groups with its usual tensor product $\otimes_\Z$.
\item Similarly, the ring $\hat{\Z}_p$ of $p$-adic integers has the structure of a unital idempotent in the category of abelian groups.
\item Any (co)localization functor (see Example \ref{ex:localization}) is by definition a (co)unital idempotent endofunctor.
\end{enumerate} 

Another important family of examples arises when categorifying structures relating to the representation theory of Hecke algebras and quantum groups.  For instance there are the categorified Jones-Wenzl idempotents \cite{CK12a,Roz10a,FSS12} and generalizations \cite{Roz10b,Rose12,Cau12,CH12,Hog15,AbHog15}.

Further, many standard constructions in homological algebra can be understood in the language of categorical idempotents.  Let $R$ be a ring, and let $\KC^-(R\modules)$ be the homotopy category of bounded above complexes of $R$-modules.  Isomorphism in $\KC^-(R\modules)$ is homotopy equivalence, and is denoted $\simeq$.  Projective resolution gives a idempotent endofunctor $\PB:\KC^-(R\modules)\rightarrow \KC^-(R\modules)$.  In fact $\PB$ has the structure of a counital idempotent in the category of triangulated endofunctors of $\KC^-(R\modules)$; the counit $\e:\PB\rightarrow \Id$ is defined by a choice of quasi-isomorphism $\e_C:\PB(C)\rightarrow C$ for each $C\in \KC^-(R\modules)$.  Localizing with respect to the quasi-isomorphisms yields the derived category $D^-(R\modules)$, which is equivalent to the the subcategory $\KC^-(R\proj)\subset \KC^-(R\modules)$ consisting of complexes of projective modules.  This category is nothing other than the essential image of $\PB$ (that is the full subcategory of complexes which are homotopy equivalent to their own projective resolutions).  To bring semi-orthogonal decompositions into the picture, standard arguments show that any complex $C\in \KC^-(R\modules)$ fits into a distinguished triangle of the form
\[
\PB(C)\rightarrow C \rightarrow \AB^-(C)\rightarrow \PB(C)[1]
\]
where $\PB(C)$ is a complex of projectives and $\AB^-(C)$ is an acyclic complex. Any such triangle is unique up to unique isomorphism, and $\AB^-(C)$ depends functorially on $C$.  Further, any chain map from a projective complex to an acyclic complex is null-homotopic.  Thus the above triangle realizes the semi-orthogonal decomposition of $\KC^-(R\modules)$ into the complexes of projectives and the acyclic complexes.  In other words, any complex can be split up into an acyclic complex and a complex of projectives in an essentially unique (and functorial) way.

Note that everything above can be described in terms of the functor $\PB$:
\begin{itemize}
\item A complex of projectives is a complex $C$ such that $\PB(C) \simeq C$.
\item A chain map $f:C\rightarrow D$ is a quasi-isomorphism iff $\PB(f)$ is a homotopy equivalence.
\item The derived category $D^-(R\modules)$ is (up to equivalence of categories) nothing other than the essential image of $\PB$; that is, the full subcategory of $\KC^-(R\modules)$ consisting of objects isomorphic to $\PB(C)$ for some $C$.
\item The functor $\AB^-$ sends a complex $C$ to the mapping cone $\Cone(\e_C:\PB(C)\rightarrow C)$.  This functor $\AB^-$ is \emph{complementary idempotent} to $\PB$, and is uniquely determined by $\PB$.
\end{itemize}

\subsection{Complementary categorical idempotents}\label{subsec:mainResults}
Let $\AS$ be a triangulated monoidal category (see \cite{Balmer10}, and also our section \S \ref{subsec:triangMonoidal}).  In particular $\AS$ has direct sum $\oplus$ and tensor product $\otimes$, suspension functor $[1]$, and a collection of \emph{distinguished triangles} $A\rightarrow B\rightarrow C\rightarrow A[1]$.  We remind the reader that the \emph{Grothendieck group} $K_0(\AS)$  of a triangulated category $\AS$ is the abelian group formally spanned by isomorphism classes  of objects of $\AS$, modulo the relation that $[A]-[B]+[C]=0$ for every distinguished triangle $A\rightarrow B\rightarrow C\rightarrow A[1]$.  If $\AS$ is a triangulated monoidal category, then $K_0(\AS)$ is a ring.  Let $\MS$ be a triangulated category on which $\AS$ acts by triangulated endofunctors.   The action will be denoted by $A(M)$ for $A\in \AS$, $M\in \MS$.  For instance we could take $\MS=\AS$ with its left regular action: $A(B):=A\otimes B$ for all $A,B\in \AS$.  The categories $\AS$ and $\MS$ will be fixed throughout.

A \emph{pair of complementary idempotents in $\AS$} is a pair of objects $(\counital,\unital)$ in $\AS$, together with a distinguished triangle
\begin{equation}\label{eq:PQtriang_intro}
\counital\buildrel\e\over \rightarrow \one \buildrel \eta\over \rightarrow \unital\buildrel \d\over \rightarrow \counital[1],
\end{equation}
subject to the condition that $\counital\otimes \unital\cong 0\cong \unital\otimes \counital$.  Note that in the Grothendieck group, the classes $[\unital]$, $[\counital]$ are a pair of complementary idempotents in $K_0(\AS)$.

In \S \ref{sec:firstProps} we establish some consequences of this definition.  For  instance, it is straightforward to prove that $\unital$ and $\counital$ have the structures of a unital and counital idempotent, respectively, and that any unital or counital idempotent has a complementary idempotent.  Further, an object absorbs an idempotent if and only if it annihilates the complement. The following is somewhat less obvious:
\begin{theorem}\label{introThm:homs}
Let $(\counital,\unital)$ be a pair of complementary idempotents and let $M,N\in \MS$ arbitrary.  We have
\begin{enumerate}
\item $\Homg_\MS(\unital(M),\unital(N))\cong \Homg_\MS(M,\unital(N))$,
\item $\Homg_\MS(\counital(M),\counital(N))\cong \Homg_\MS(\counital(M),N)$, and
\item $\Homg_\MS(\counital(M),\unital(N))\cong 0$.
\end{enumerate}
In fact the isomorphism (1) is induced by precomposing with the map $\one\rightarrow \unital$, and the isomorphism (2) is induced by post-composing with $\counital\rightarrow \one$.  
\end{theorem}
Here, $\Homg_\AS(A,B)$ denotes the graded space of homs: $\Homg_\AS(A,B):=\bigoplus_{i\in\Z} \Homg_\AS^i(A,B)$, where $\Homg_\AS^i(A,B):=\Hom_\AS(A,B[i])$.  A particularly important special case occurs when $\MS=\AS$, and $M=N=\one$:
\begin{enumerate}
\item[(1')] $\Homg_\AS(\unital,\unital)\cong \Homg_\AS(\one,\unital)$.
\item[(2')] $\Homg_\AS(\counital,\counital)\cong \Homg_\AS(\counital,\one)$.
\item[(3')] $\Homg_\AS(\counital,\unital)\cong 0$.
\end{enumerate}
Each of the isomorphisms (1) and (2) is an isomorphism of graded commutative algebras.    This theorem is restated and proven in \S \ref{subsec-semiOrtho}.  

Motivated by the example of group cohomology, we make the following:

\begin{definition}\label{def:cohomology}
Let $(\counital,\unital)$ be a pair of complementary idempotents in $\AS$.  The cohomology of $\unital$ (resp.~$\counital$) is the graded commutative algebra from statement (1') (resp.~(2')) above.  If $\MS$ is a category on which $\AS$ acts by triangulated endofunctors, then the \emph{cohomology of $\unital$ and $\counital$ with coeffiecients in $M,N\in \MS$} are the graded spaces from statements (1) and (2) of Theorem \ref{introThm:homs}, respectively. 
\end{definition}

\begin{example}\label{ex:introHopf}
Let $\k$ be a commutative ring and let $H$ be a Hopf algebra (e.g.~$H=\k[G]$ the group algebra of a discrete group).  The category of $H$-modules is monoidal with tensor product $\otimes=\otimes_\k$ and monoidal identity given by the trivial module $\one =\k$.  Let $\e:\PB\rightarrow \k$ be a resolution of the trivial module by projective $H$-modules.  Then $(\PB,\e)$ is a counital idempotent in $\KC^-(H\modules)$.  The cohomology of $\PB$ is isomorphic to $\Ext_H(\k,\k)$, which by definition is the usual cohomology of $H$ (if $H=\k[G]$, then this is also called the group cohomology of $G$), by definition.
\end{example}

\begin{example}\label{ex:introBimod}
Let $\k$ be a commutative ring, let $R$ be a $\k$-algebra which is flat as a $\k$-module, and set $R^e:=R\otimes_\k R^{\op}$.  An $(R,R)$-bimodule is the same thing as an $R^e$-module, hence the category of $R^e$-modules is monoidal with tensor product $\otimes =\otimes_R$ and monoidal identity given by the trivial bimodule $\one =R$.  Let ${\PB}\rightarrow R$ be a resolution of $R$ by projective $R^e$-modules.  Then ${\PB}$ is a counital idempotent in $\KC^-(R^e\modules)$.  The cohomology of ${\PB}$ is isomorphic to $\Ext_{R^e}(R,R)$, which is by definition the Hochschild cohomology of $R$.
\end{example}

We refer to statement (3) of Theorem \ref{introThm:homs} as the \emph{semi-orthogonality of categorical idempotents}.  The semi-orthogonality property forces a certain rigidity for categorical idempotents which does not exist for ordinary idempotents in linear algebra.  To illustrate, note that in linear algebra the image of an idempotent does not determine its kernel, and vice versa.  For example, there are infinitely many idempotent operators $\R^2\rightarrow \R^2$ whose image is a given line $\R\cong L\subset \R^2$.  Note, however, that there is a unique \emph{orthogonal projection} onto a given line.  Categorical idempotents have a sort of orthogonality built in, one consequence of which is the following:

\begin{theorem}
Let $\idempotent$ be a unital or counital idempotent in $\AS$.  Then $\idempotent$ is characterized up to isomorphism by either of the following full subcategories of $\AS$:
\begin{enumerate}
\item The full subcategory consisting of objects $A\in \AS$, such that $\idempotent\otimes A \cong A$.
\item The full subcategory consisting of objects $A\in \AS$, such that $\idempotent\otimes A \cong 0$.
\item The full subcategory consisting of objects $A\in \AS$, such that $A\otimes \idempotent \cong A$.
\item The full subcategory consisting of objects $A\in \AS$, such that $A\otimes \idempotent \cong 0$.
\end{enumerate}
In fact, if $(\unital,\counital)$ is a pair of complementary idempotents in $\AS$,  then the triangle
\[
\counital \rightarrow \one\rightarrow \unital \rightarrow \counital[1]
\]
is uniquely determined by any of the above categories, up to unique isomorphism of triangles.
\end{theorem}
This theorem is a consequence of the fundamental Theorem \ref{thm:idemptOrder}, which we restate in a weaker form here:
\begin{theorem}\label{introThm:partialOrder}
Let $(\unital_i,\eta_i)$ be unital idempotents ($i=1,2$).  The following are equivalent:
\begin{enumerate}
\item $\unital_1\otimes \unital_2\simeq \unital_1$.
\item $\unital_2\otimes \unital_1\simeq \unital_1$.
\item There exists a map $\nu:\unital_1\rightarrow \unital_2$ such that $\nu\circ \eta_1 = \eta_2$.
\end{enumerate}
Furthermore, if any of these equivalent statements is true, then the map $\nu$ is unique, and there is a unique morphism of triangles
\[
\begin{diagram}[small]
\counital_1 & \rTo^{\e_1} & \one & \rTo^{\eta_1} & \unital_1 & \rTo^{\d_1} & \counital_1[1]\\
\dTo && \dTo^{\Id} && \dTo^{\nu} && \dTo \\ 
\counital_2 & \rTo^{\e_2} & \one & \rTo^{\eta_2} & \unital_2 & \rTo^{\d_2} & \counital_2[1]\\
\end{diagram}\]
Similar statements hold for counital idempotents.
\end{theorem}
This is restated and proven in \S \ref{subsec-uniqunessOfIdempts}.  This allows us to put a partial order on the set of isomorphism classes of unital idempotents (following \cite{BD14}) by $\unital_1\leq \unital_2$ if $\unital_1\otimes \unital_2\simeq \unital_1$.

\begin{remark}\label{rmk:relativeIdempotents}
We refer to statement (3) of Theorem \ref{introThm:partialOrder} as the fact that $\unital_2$ is a \emph{unital idempotent relative to $\unital_1$}.  This special relationship implies, for instance, that $\Homg_\AS(\unital_1,\unital_2)\cong \Homg_\AS(\unital_2,\unital_2)$; compare with Theorem \ref{introThm:homs}.
\end{remark}

\subsection{Idempotents and semi-orthogonal decompositions}
\label{subsec:introDeomp}
In this section we explain what is meant by sem-orthogonal decomopsition, and what is its relationship with categorified idempotents.  This is a valuable perspective

Let $\AC,\BC$ be categories and let $F:\AC\rightarrow \BC$ be a functor with a right adjoint $G:\BC\rightarrow \AC$.  It is well known that $G$ is full if and only if the unit of the adjunction $\eta:\Id\rightarrow G\circ F$ makes $G\circ F$ into a unital idempotent in $\End(\BC)$.  Similarly, $F$ is full if and only if the counit of the adjuntion $\e:F\circ G\rightarrow \Id$ makes $F\circ G$ into a counital idempotent in $\End(\AC)$.  This is best understood in the language of monads (see Proposition 4.2.3 in \cite{BorceuxBook2}). 

Now, suppose that $\MS$ is a triangulated category, and let $0=\NS_0\subset \NS_1\subset \cdots \NS_r = \MS$ be a sequence of full, triangulated subcategories.  The sequence of subcategories $0=\NS_0\subset \NS_1\subset \cdots \NS_r = \MS$ is called a \emph{semi-orthogonal decomposition} of $\MS$ if the inclusion $\IB_i:\NS_i\rightarrow \MS$ admits a left adjoint $\PB_i:\MS\rightarrow \NS_i$.  Suppose this is so, and set $\UB_i:=\IB_i\circ \PB_i$.   As in the above discussion, the unit of the adjunction makes $\unital_i$ into a unital idempotent in the monoidal category of exact endofunctors $\End(\MS)$.  It is trivial to check that if $i\leq j$, then $\unital_j\leq \unital_i$.

\begin{example}[Bousfield localization]\label{ex:localization}
Let $\AS$ be a triangulated monoidal category which acts on a triangulated category $\MS$ by triangulated endofunctors.  If $(\unital,\eta)$ is a unital idempotent in $\AS$, then the action of $\unital$ on $\MS$ is called a \emph{localization functor}.   An object $Z\in \MS$ is \emph{$\unital$-acyclic} if $\unital(Z)\cong 0$.  A morphism $f$ in $\MS$ is called a $\unital$-quasi-isomorphism if $\unital(f)$ is an isomorphism.  Equivalently, $f$ is a $\unital$-quasi-isomorphism  iff $\Cone(f)$ is $\unital$-acyclic.  Localizing with respect to the $\unital$-quasi-isomorphisms yields the \emph{Bousfield localization} of $\MS$ with respect to $\unital$.  The resulting category is equivalent to the essential image $\im\unital\subset \MS$.  We denote this image by $\unital\MS$; it is the category consisting of the objects $X\in \MS$ which are \emph{$\unital$-local}, i.e. $\unital(X)\cong X$.

The complementary idempotent $\unital^c$ gives rise to a \emph{colocalization functor} $\MS\rightarrow \MS$.  For each $X\in \MS$ there is a distinguished triangle
\begin{equation}\label{eq:locTriang}
\unital^c(X)\rightarrow X \rightarrow \unital(X)\rightarrow \unital^c(X)[1].
\end{equation}
These distinguished triangles are functorial in $X$.   Indeed, by definition there is a distinguished triangle already in $\AS$:
\[
\unital^c \rightarrow \one \rightarrow \unital\rightarrow \unital^c[1].
\]
Evaluating on $X$ gives the distinguished triangle (\ref{eq:locTriang}).  If $Y\in \unital\MS$, then the isomorphism $\Hom_{\unital\MS}(\unital X,Y)\cong \Hom_{\MS}(X,Y)$ from Theorem \ref{introThm:homs} implies that the inclusion $\unital\MS\rightarrow \MS$ has a left adjoint $\MS\rightarrow \unital\MS$ sending $X\mapsto \unital X$.  Thus, Bousfield localization is equivalent to providing a short semi-orthogonal decomposition $0=\NS_0\subset \NS_1\subset \NS_2=\MS$, where $\NS_1=\unital\MS_1$.

The phrase semi-orthogonal is justified by Theorem \ref{introThm:homs}, which states that there are no nonzero homs from a $\unital$-acyclic object $X\in \im \unital^c = \ker \unital$ to a $\unital$-local object $Y\in \im \unital$. 
\end{example}

\begin{example}
Let $\AS$ be a triangulated monoidal category which acts on a triangulated category $\MS$ by triangulated endofunctors. Let $\{(\unital_i,\eta_i)\}_{i=0}^r$ be a sequence of unital idempotents with $0=\unital_0\leq \unital_1\leq \cdots \leq \unital_r\cong \one$.  Then $\NS_i:=\unital_i\MS$ gives a semi-orthogonal decomposition of $\MS$. 
\end{example}

Motivated by this example, we define an idempotent decomposition of $\one$ in a triangulated monoidal category to be a particular sort of partially ordered family of unital idempotents.  We explore consequences of this definition in \S \ref{subsec:generalizedDecomp}.

In these notes we focus on studying the idempotents in $\AS$ themselves, rather than how they act on $\MS$.  Thus, we will not mention semi-orthogonal decomposition or localization further.  

\subsection{Tate cohomology}
\label{subsec:introTate}
We were originally led to consider Tate cohomology by the following.

\begin{question}\label{question:HomAlg}
Is there a natural algebra structure on $\Homg_\AS(\unital,\counital)$?
\end{question}

One naive attempt at defining such an algebra structure would be the following: given $f,g\in \Homg_\AS(\unital,\counital)$, define their product to be the composition
\[
\begin{diagram}[small]
f\ast g\ \ \:=  \ \ \ \unital &\rTo^{\cong} & \unital\otimes \unital & \rTo^{f\otimes g} & \counital\otimes \counital &\rTo^{\cong} & \counital
\end{diagram}
\]
However under this definition $f\ast g=0$ for all $f,g\in \Homg_\AS(\unital,\counital)$, since $f\ast g$ factors through $\unital\otimes \counital$, which is isomorphic to zero.  As it turns out, under certain conditions $\Homg_\AS(\unital,\counital[1])$ is an algebra with unit given by the connecting morphsim $\d:\unital\rightarrow \counital[1]$ (see Theorem \ref{thm:TateAndDuals} below).  We will arrive at this result somewhat indirectly.

Let $\e:\counital\rightarrow \one$ be a counital idempotent.  Suppose we are given a \emph{unital} idempotent $\e^\star:\one\rightarrow \counital^\star$ which is to be compared with $\counital$.  A priori there is no relation assumed between $\counital$ and $\counital^\star$, though the notation is suggestive for later purposes.  To study the relationship between $\counital$ and $\counital^\star$, it is convenient to introduce an object $\Tate(\counital,\counital^\star)$ which is defined to be the mapping cone on the map $\e^\star\circ \e:\counital\rightarrow \counital^\star$.  Similarly, one has an object $\Tate(\unital^\star,\unital)$ where $\unital$ and $\unital^\star$ are the idempotents complementary to $\counital$ and $\counital^\star$.  In \S \ref{subsec:generalTate} we prove the following:

\begin{theorem}\label{introThm:Tate}
We have $\Tate(\counital,\counital^\star)\simeq \Tate(\unital^\star,\unital)$.  Henceforth we will denote this object simply by $\Tate$.  Furthermore:
\begin{enumerate}
\item  If $\counital\otimes \unital^\star\simeq 0\simeq \unital^\star\otimes \counital$, then $\Tate\simeq \unital\otimes \counital^\star$ has the structure of a unital idempotent in $\AS$, with complementary idempotent $\counital\oplus \unital^\star$.
\item If $\unital\otimes \counital^\star\simeq 0 \simeq \counital^\star\otimes \unital$, then $\Tate[-1]\simeq \counital\otimes\unital^\star$ has the structure of a counital idempotent in $\AS$ with complementary idempotent $\unital\oplus \counital^\star$.
\end{enumerate}
\end{theorem}

We call $\Tate$ the Tate object associated to $(\counital,\counital^\star)$ (or, equivalently, assocated to $(\unital^\star,\unital)$).  

\begin{definition}\label{def:TateCohomology}
The \emph{Tate cohomology of the pair $(\counital,\counital^\star)$} (or, equivalently of the pair $(\unital^\star,\unital)$) is the graded algebra $\Homg_\AS(\Tate,\Tate)$.  If $\AS$ acts on $\MS$, then the \emph{Tate cohomology of $(\counital,\counital^\star)$ with coefficients in $M,N\in \MS$} is the graded space $\Homg_\MS(\Tate(M),\Tate(N))$.
\end{definition}
Of course, if we are in the situation (1) or (2) of Theorem \ref{introThm:Tate}, then the Tate cohomology can be simplified according to Theorem \ref{introThm:homs}.

\begin{example}\label{ex:introTate}
We summarize the usual set-up for Tate cohomology, referring to the classic preprint of Buchweitz \cite{BuchweitzPreprint} for the details.  Let $S$ be a strongly Gorenstein ring and $\MS=\KC^f(S)$ the homotopy category of unbounded complexes of  (right) $S$-modules, whose homology is finitely generated as an ungraded $S$-module.  Projective resolution determines a functor $\PB:\MS\rightarrow \MS$ (for projective resolution of unbounded complexes, see \cite{Spalt_resolutions}).  There is also a projective \emph{coresolution} functor $\CB:\MS\rightarrow \MS$.  These functors are equipped with natural transformations $\e:\PB\rightarrow \Id_{\MS}$ and $\eta:\Id_{\MS}\rightarrow \CB$ such that $(\PB,\e)$ and $(\CB,\eta)$ are a counital and a unital idempotent in $\End(\MS)$.

The Tate object associated to $\PB,\CB$ is a shift of the \emph{complete resolution} functor:
\[
\Tate(M):=\Cone(\PB(M)\rightarrow \CB(M))\cong \CB\RB(M)[1]
\]
where the map $\PB(M)\rightarrow \CB(M)$ the evident composition of quasi-isomorphisms.  The homology of the complex $\Hom_S^\bullet(\CB\RB(M),\CB\RB(N))$ is the Tate cohomology of $S$ with coefficients in $M,N$.  In the language of this paper, it could also be referred to as the Tate cohomology of $(\PB,\CB)$ with coefficients in $M,N\in \KC^f(S)$.  

 Let $\AB^-(M):=\Cone(\PB(M)\rightarrow M)$.  Then $\AB^-$ is the complementary idempotent to $\PB$.  It is possible to show that $\AB^-\circ \CB\cong 0 \cong \CB\circ \AB^-$, so that statement (2) of Theorem \ref{introThm:Tate} applies.  Thus, complete resolution $\CB\RB\cong \Tate[-1]$ is a counital idempotent endofunctor, and Tate cohomology with coefficients in $M,N$ can also be computed as the homology of the complex of homs $\Hom_S^\bullet(\CB\RB(M),N)$.  The essential image of $\CB\RB$ is also called the stable category of $S$.
\end{example}

Now we specialize to the case where $\AS$ is a homotopy category of complexes.  Let $\CS$ be an additive monoidal category.  In general the tensor product of unbounded complexes in $\KC(\CS)$ is undefined.  However, $\KC(\CS)$ includes as a full subcategory into the monoidal categories $\KC(\CS^\oplus)$ and $\KC(\CS^\Pi)$, where $\CS^\oplus$ and $\CS^\Pi$ are the categories obtained from $\CS$ by formally adjoining countable direct sums and products, respectively.  The tensor products in these categories are
\begin{equation}\label{eq:tensorProdSum}
(A\otimes B)_k =\bigoplus_{i+j=k} A_i\otimes B_j \ \ \  \ \ \text{ with differential } \ \ \ \ \ \ d|_{A_i\otimes B_j} = d_{A_i}\otimes \Id_{B_i} + (-1)^{i}\Id_{A_i}\otimes d_{B_j}
\end{equation}
and
\begin{equation}\label{eq:tensorProdProd}
(A\cotimes B)_k =\prod_{i+j=k} A_i\otimes B_j \ \ \  \ \ \text{ with differential } \ \ \ \ \ \ d|_{A_i\otimes B_j} = d_{A_i}\otimes \Id_{B_i} + (-1)^{i}\Id_{A_i}\otimes d_{B_j},
\end{equation}
respectively.  Thus, there are \emph{two} ways of tensoring unbounded complexes: one using direct sums, and the other using direct products.

\begin{theorem}\label{thm:TateAndDuals}
Let $\CS$ be a rigid, additive monoidal category with duality $(-)^\star$.  Let $(\counital, \unital)$ be a pair of complementary idempotents in $\KC^-(\CS)$, and let
\[
\Tate=\Cone(\counital\rightarrow \counital^\star)\simeq \Cone(\unital^\star\rightarrow \unital) \ \ \ \ \in \ \ \ \  \KC(\CS)
\]
be the Tate object.  Then:
\begin{enumerate}
\item Inside $\KC(\CS^{\oplus})$, we have $\counital\otimes \unital^\star\simeq 0 \simeq \unital^\star\otimes \counital$, hence $\Tate$ is a unital idempotent in $\KC(\CS^\oplus)$.
\item Inside $\KC(\CS^\Pi)$, we have $\unital\cotimes \counital^\star\simeq 0 \simeq \counital^\star\cotimes \unital$, hence $\Tate[-1]$ is a counital idempotent in $\KC(\CS^\Pi)$.
\end{enumerate}
Furthermore, $\Homg_\AS(\unital,\counital)\cong \Homg_\AS(\counital^\star,\unital^\star)\cong \Homg_\AS(\Tate,\Tate)$; each is a graded commutative algebra with unit given by class of the connecting morphism $\unital\rightarrow \counital[1]$, respectively $\counital^\star\rightarrow \unital^\star[1]$, respectively $\Id_\Tate$.
\end{theorem}
  The adjective ``rigid'' means that there is a duality functor $(-)^\star:\CS\rightarrow \CS^{\op}$ such that, among other things, $\Hom_\CS(a,b)\cong\Hom_\CS(a\otimes b^\star, \one)\cong \Hom_\CS(\one,a^\star\otimes b)$.  Extending $(-)^\star$ to complexes gives a pair of inverse contravariant functors $\KC^\pm(\CS)\rightarrow \KC^\mp(\CS)$ and $\KC(\CS^\oplus)\leftrightarrow \KC(\CS^\Pi)$.




\subsection{Acknowledgements}
As mentioned already, most of this material is presumably well-known in certain circles.  Our intention with these notes is to make it available to a broader audience, and not necessarily to claim ownership of any particular result.   That said, we were particularly influenced by \cite{BD14} at various points, especially as regards the partial order on idempotents in \S \ref{subsec-uniqunessOfIdempts}.  Our perspective on Tate cohomology is similar to that of Greenlees \cite{Greenlees87_tate} and Greenlees-May \cite{GreenleesMay95}, which we learned of after most of this work was written.

\subsection{Notation}
If $\CC$ is an additive category, then we let $\KC(\CC)$ denote the homotopy category of unbounded complexes over $\CC$ (differentials will have degree $+1$).  Let $\KC^{-,+,b}(\CC)\subset \KC(\CC)$ denote the full subcategories of bounded above, bounded below, respectively bounded complexes.  If $\CC$ is a category with an invertible endofunctor (for instance, grading shift) $[1]$, we set $\Homg_\CC(A,B):=\bigoplus_k \Hom_\CC(A,B[k])$.  In the special case of complexes, $[1]$ will denote the downward grading shift: $A[1]_i=A_{i+1}$, with differential $d_{A[1]}=-d_A$.  The sign is conventional, and mostly irrelevant since $(C,d_C)\cong (C,-d_C)$ for any chain complex $C$ with differential $d_C$.  If $f:A\rightarrow B$ is a chain map, then the mapping cone is $\Cone(f)=A[1]\oplus B$ with differential $\smMatrix{-d_A & 0 \\ f & d_B}$.  The sign convention for $[1]$ ensures that the projection map $\Cone(f)\rightarrow A[1]$ is a chain map.

\section{A motivating example}
\label{sec:mainExample}
In this section we introduce what we feel is the simplest nontrivial example of a pair of complementary idempotents.   All of the essential features of the theory are already present in this simple case.

Let $\k$ be a commutative ring and set $C:=\k[\Z/2]=k[x]/(x^2-1)$.  We regard $\k$ as a $C$-module, where $x$ acts by 1.  Let $\AC$ be the category of left $C$-modules.  Given $M,N\in \AC$, the tensor product $M\otimes_\k N$ inherits a $C$-action via $x\mapsto x\otimes x$.   Consider the following complex of $C$-modules:
\[
\PB \ \ = \ \ \begin{diagram}[small] \cdots & \rTo^{1+x} &C& \rTo^{1-x} &C& \rTo^{1+x} &C & \rTo^{1-x} & \underline{C},\end{diagram},
\]
in which we have underlined the term in homological degree zero.  Observe that there is a chain map $\PB\rightarrow \k$ which describes a projective resolution of the trivial $C$-module $\k$.  By uniqueness of projective resolutions, we have $\PB\otimes_\k\PB\simeq \PB$, so $\PB$ is an idempotent in $\KC^-(\AC)$.  It is natural ask if there is an idempotent which complementary to $\PB$.  In fact there is.  Consider the complex
\[
\AB \ \ = \ \ \begin{diagram}[small] \cdots & \rTo^{1-x} &C& \rTo^{1+x} &C& \rTo^{1-x} &C & \rTo^{} & \underline{\k},\end{diagram}
\]
where the map $C\rightarrow \k$ is the unique $C$-module map sending $x\mapsto 1$.  Observe that $\AB$ is acyclic, being the cone on the quasi-isomorphism $\PB\rightarrow \k$.  We emphasize, however, that $\AB$ is nontrivial in $\KC^-(\AC)$, since acyclic complexes need not be contractible (exact sequences of $C$-modules need not be split exact).  There is a short exact sequence of complexes $0\rightarrow \k \rightarrow \AB\rightarrow \PB[1]\rightarrow 0$, which gives rise to a distinguished triangle
\[
\PB \rightarrow \k \rightarrow \AB \rightarrow \PB[1]
\]
in $\KC^-(\AC)$.  This expresses the fact that $\AB$ is \emph{complementary} to $\PB$.  In order to prove that $\AB\otimes_\k\AB\simeq \AB$, it will suffice to prove that $\PB\otimes_\k \AB\simeq 0\simeq \AB\otimes_\k\PB$ (see Proposition \ref{prop-projectorAbsorbing}).  To see this, observe that $C\otimes_\k M$ is a projective $C$-module for any $M$, hence $\PB\otimes_\k\AB$ is a complex of projective $C$-modules.  On the other hand $C$ is free as a $\k$-module, so $C\otimes_\k \AB$ is acyclic.  Thus, $\PB\otimes_\k \AB$ is a bounded above, acyclic complex of projectives.  It follows that $\PB\otimes_\k\AB\simeq 0$; this is simply a restatement of the usual fact that a semi-infinite exact sequence $\cdots \rightarrow P_{-1}\rightarrow P_0\rightarrow 0$ of projectives is split exact.  In Definition \ref{def:idempotents} Thus, $(\PB,\AB)$ is a pair of complementary idempotents in the sense of Definition \ref{def:idempotents}.

Retain notation as above.  If $M$ is a $C$-module, then the linear dual $M^\ast:=\Hom_\k(M,\k)$ is a $C$-module via
\[
x\phi \ : \ m\mapsto \phi(xm)
\]
for all $\phi\in M^\ast$.  Recall the complexes  $\PB$ and $\AB$.  Their duals are the complexes
\[
\PB^\ast  \ \ = \ \ \begin{diagram}[small] \underline{C} &\rTo^{1-x}& C &\rTo^{1+x}& C &\rTo^{1-x}& C &\rTo^{1+x}& \cdots\end{diagram},
\]
respectively
\[
\AB^\ast \ \ = \ \ \begin{diagram}[small]  \underline{\k} &\rTo^{}& C &\rTo^{1-x}& C &\rTo^{1+x}& C &\rTo^{1-x}&\cdots \end{diagram} .
\]

It is not hard to show that
\begin{proposition}\label{prop:specialExample}
The following relations hold:
\begin{enumerate}
\item $\AB\otimes \PB\simeq 0 \simeq \PB\otimes \AB$,
\item $\AB^\ast\otimes \PB^\ast \simeq 0 \simeq \PB^\ast\otimes \AB^\ast$,
\item $\AB^\ast \otimes \PB \simeq 0 \simeq \PB\otimes \AB^\ast$.
\end{enumerate}
\end{proposition}
Note that this last tensor product requires the formation of countable direct sums.  Statement (1) was proven already, and statement (2) follows from (1) by duality.  We leave the proof of statement (3) as an exercise, since it follows from our work in \S \ref{subsec:duals}.  
It is quite conspicuous that we have up to this point said nothing about the tensor product $\AB\otimes \PB^\ast\cong\PB^\ast\otimes \AB$.  It will turn out that this complex is quite interesting, and not contractible.

Note that there is a chain map $\PB\rightarrow \PB^\ast$.  We will denote the mapping cone on this map on this map by $\Tate$.  Explicitly,
\[
\Tate \ \ :=  \ \  \begin{diagram}[small] \cdots &\rTo^{1+x} & C &\rTo^{1-x}& \underline{C} &\rTo^{1+x}& \underline{C} &\rTo^{1-x}& C &\rTo^{1+x}& \cdots\end{diagram}.
\]
Note that $\Tate$ is also equivalent to the mapping cone on a map $\AB^\ast\rightarrow \AB$.
\[
\Tate \ \ \simeq  \ \  \begin{diagram}[small]
 & & & & \k&\rTo & C&\rTo & C&\rTo & \cdots\\
 & & & & \oplus &\rdTo^{\Id} & \oplus & & & & \\
\cdots &\rTo & C&\rTo & C&\rTo & \k& & & & 
\end{diagram}.
\]
\begin{proposition}\label{prop:tateIdempt1}
We have $\AB\otimes \PB^\ast\simeq \Tate\simeq \PB^\ast\otimes \AB$.
\end{proposition}
\begin{proof}
Consider the exact triangle
\[
\PB\rightarrow \k \rightarrow \AB\rightarrow \PB[1].
\]
Tensoring with $\AB^\ast$ and using the fact that $\AB^\ast\otimes \PB\simeq 0$ gives a distinguished triangle
\[
0\rightarrow \AB^\ast\rightarrow \AB^\ast\otimes \AB\rightarrow 0.
\]
From which it follows that $\AB^\ast\otimes \AB\simeq \AB^\ast$.  A similar argument shows that $\PB^\ast\otimes \PB\simeq \PB$. 

Now, $\AB$ can be expressed as the mapping cone on the map $\PB\rightarrow \k$.  Tensoring with $\PB^\ast$, we see that $\PB^\ast\otimes A$ is homotopy equivalent to the mapping cone on a map $\PB\rightarrow \PB^\ast$.  Careful examination reveals that this map is the composition of canonical maps $\PB\rightarrow \k\rightarrow \PB^\ast$, and the proposition now follows by definition of $\Tate$.
\end{proof}

We leave it as an exercise that there is some differential on $\Tate\oplus \AB^\ast\oplus \PB$ so that the resulting complex is homotopy equivalent to $\k$.  In fact, the only nonzero components of the differential (besides those differentials internal to each summand) go from $\Tate$ to $\AB^\ast$ or $\PB$.  This information can be represented somewhat schematically by
\[
\k \ \simeq \ \Tot\left(
\begin{tikzpicture}[baseline=0cm]
\node(a) at (0,0) {$\Tate$};
\node(b) at (2,1) {$\AB^\ast$};
\node(c) at (2,-1) {$\PB$};
\tikzstyle{every node}=[font=\small]
\path[->,>=stealth',shorten >=1pt,auto,thick,node distance=1.8cm]
(a) edge node[above] {$[1]$} (b)
(a) edge node[above] {$[1]$} (c);
\end{tikzpicture}
\right).
\]
The three terms above are mutually orthogonal, hence we obtain an interesting decomposition of the identity module $\k$ in to three mutually orthogonal idempotents with respect to $\otimes$.

There is an alternate definition of the tensor product (see Equation (\ref{eq:tensorProdProd}) in which we replace direct sum by direct product.  We denote this alternate tensor product by $\cotimes$.  If $M$ and $N$ are bounded above (or bounded below), then $M\otimes N\cong M\cotimes N$.  It is an exercise to show that using $\cotimes$ we have $\AB\cotimes \PB^\ast\simeq 0 \simeq \PB^\ast\cotimes \AB$, then to verify that $\AB^\ast\cotimes \PB\simeq \Tate[-1]\simeq \PB\cotimes \AB^\ast$.  Furthermore, 
\[
\k \ \simeq \ \left(
\begin{tikzpicture}[baseline=0cm]
\node(a) at (0,0) {$\Tate[-1]$};
\node(b) at (-2,1) {$\AB$};
\node(c) at (-2,-1) {$\PB^\ast$};
\tikzstyle{every node}=[font=\small]
\path[->,>=stealth',shorten >=1pt,auto,thick,node distance=1.8cm]
(b) edge node[above] {$[1]$} (a)
(c) edge node[above] {$[1]$} (a);
\end{tikzpicture}
\right)
\]
is a decomposition of $\k$ into three mutually orthogonal idempotents (with respect to the ``direct product'' version $\cotimes$ of the tensor product).

\section{Triangulated monoidal categories}
\label{sec:triangMonoidal}

Recall that if $e\in A$ is an idempotent element of a unital algebra, then the complementary idempotent is $1-e\in A$.  In order to have a categorical analogue of complementary idempotent we will therefore pass to the setting of \emph{triangulated monoidal categories}, where the role subtraction is played by the formation of mapping cones.  Thus, we begin with a brief reminder on the basics of triangulated categories, monoidal categories, and how they will interact.

\subsection{Graded monoidal categories}
\label{subsec:gradedMonoidal}
We assume the reader is comfortable with the basics of monoidal categories (see for instance Chapter 2 of \cite{EGNOtensorCat}).  A monoidal category is a 6-tuple $(\AS,\otimes,\one,\mu,\l,\rho)$ where $\AS$ is a category $\otimes:\AS\times \AS\rightarrow \AS$ is a functor, $\one\in \AS$ is an object, and $\mu,\l,\rho$ are certain natural isomorphisms (called the associator and left/right unitors).  These are subject to certain conditions called coherence relations.

There are certain signs that appear in the compatibility relations between categories with triangulated and monoidal structures.  We give an explanation of these signs first in an easier context: that of graded monoidal categories.  

\begin{definition}\label{def:gradedMonoidalCat}
A \emph{graded monoidal category} is a 9-tuple $(\AS,\otimes,\one,\mu,\l,\rho, [1],\a,\b)$ where $\AS$ is an additive category, $(\AS,\otimes,\one,\mu,\l,\rho)$ is a monoidal category, $[1]:\AS\rightarrow \AS$ is an invertible endofunctor called the suspension (or grading shift), and:
\begin{enumerate}
\item $\a$ is a natural isomorphism $\a_A:A\otimes \one[1] \buildrel \cong\over\rightarrow A[1]$, which we call the \emph{shift absorption isomorphism}, and
\item $\b$ is a natural isomorphism $\b_A: \one[1]\otimes A\buildrel\cong\over\rightarrow A\otimes \one[1]$, which we call the \emph{braiding morphism}.
\end{enumerate}
These data are subject to the following constraints (omitting parentheses and associators for readability):

\begin{itemize}
\item[(GM1)] $\a_\one = \l_{\one[1]}:\one\otimes \one[1] \rightarrow \one[1]$, and the following diagram commutes:
\begin{equation}\label{eq:absorbtor}
\begin{diagram}
A\otimes \one[1]\otimes \one[1] & \rTo^{\a_A\otimes \Id_{\one[1]}} &A[1]\otimes  \one[1]\\
\dTo^{\a_{A\otimes \one[1]}} && \dTo_{\a_{A[1]}}\\
(A\otimes \one[1])[1] & \rTo^{\a_A[1]} & A[2]
\end{diagram}
\end{equation}
\item[(GM2)]$\b_{\one[1]} = \pi\Id_{\one[1]\otimes\one[1]}$ where $\pi=\pm 1$, and the following diagram commutes:
\begin{equation}\label{eq:commutor}
\begin{diagram}
\one[1]\otimes A\otimes B   & \rTo^{\b_A\otimes \Id_B} & A\otimes \one[1] \otimes B \\
& \rdTo^{\b_{A\otimes B}} & \dTo_{\Id_A\otimes \b_B}\\
&& A\otimes B \otimes \one[1]
\end{diagram},
\end{equation}
\end{itemize}
The scalar $\pi$ may be called the \emph{signature} of the graded monoidal category.  If $\pi=-1$, we will say that the graded monoidal category is \emph{of homological type}.  By abuse, we often simply refer to $\AS$ as a graded monoidal category.
\end{definition}

\begin{remark}
If $\AS$ is $\k$-linear for some commutative ring $\k$, then in (GM2) we may allow $\pi\in \k$ to be any scalar with $\pi^2=1$.  Unless specified otherwise, all graded monoidal categories will be of homological type ($\pi=-1$).
\end{remark}

\begin{example}\label{ex:sign}
Let $\k$ be a commutative ring, and let $\k\modules$ be the category of $\k$-modules.  Then $\k\modules$ has the structure of a monoidal category with monoid $\otimes=\otimes_\k$ and monoidal identity $\one=\k$.  The homotopy categories $\KC^{+,-,b}(\k\modules)$ each inherit the structures of graded monoidal categories from $\k\modules$, with with $[1]$ given by the usual suspension of complexes: $A[1]_i=A_{i+1}$.  For any complex $A\in \KC(\k\modules)$, the sign rule is such that
\[
d_{\k[1]\otimes A }(1\otimes a) = -1\otimes d(a) \hskip1.2in d_{A\otimes \k[1]}(a\otimes 1) = d(a)\otimes 1
\]
for all $a\in A$.  Thus, the braiding morphism $\b_A:\k[1]\otimes A \rightarrow A\otimes \k[1]$ must involve a sign.  By convention, we choose $\b_A(a\otimes 1) := (-1)^{|a|} 1\otimes a$.  Specializing to $A=\k[1]$ we get the sign in (GM2).
\end{example}

Using the associator, the absorption isomorphism, and the braiding morphism, we obtain natural isomorphisms
\[
A[k]\cong A\otimes\one[k] \ \ \ \ \ \ \ \ \ \one[k]\otimes A\cong A\otimes \one[k] \ \ \ \ \ \ \ \ A[k]\otimes B\cong (A\otimes B)[k]\cong A\otimes B[k].
\]

\begin{definition}\label{def:enrichedHom}
Let $\AS$ be a graded monoidal category, and $A,B\in \AS$ arbitrary.  Define $\Homg_{\AS}(A,B):=\bigoplus_{k\in \Z}\Hom_{\AS}(A,B[k])$ and $\Endg_\AS(A):=\Homg_\AS(A,A)$.  If $f\in \Hom_{\AS}(A,B[k])$, then we will think of $f$ as a morphism from $A$ to $B$ of degree $k$, and we write $|f|=k$.  The composition of $f:A\rightarrow B[k]$ and $g:B\rightarrow C[\ell]$ is by definition $g[k]\circ f : A\rightarrow C[k+\ell]$.  In this way, one may consider the category $\AS^\bullet$ with the same objects as $\AS$, but with hom spaces given by the \emph{enriched hom space} $\Homg_{\AS}(A,B)$.  The tensor product lifts to a bilinear map
\[
\Homg_\AS(A,A')\times \Homg_\AS(B,B') \rightarrow \Homg_\AS(A\otimes A',B\otimes B')
\]
which sends $(f,g)$ to the composition
\[
\begin{diagram}[small]
A\otimes B &\rTo^{f\otimes g} & A'[k]\otimes B'[\ell] & \rTo^{\cong}& (A'\otimes B')[k+\ell].
\end{diagram}
\]
where $k=|f|$ and $\ell=|g|$.
\end{definition}

\begin{proposition}\label{prop:superCommutativity}
Let $\AC$ be a graded monoidal category.  Then the tensor product and composition of enriched homs satisfies:
\begin{equation}\label{eq:signInTensorComp}
(f'\otimes g')\circ (f\otimes g) = (-1)^{|g'||f|}(f'\circ f)\otimes (g'\circ g).
\end{equation}
\end{proposition}
\begin{proof}
This is straightforward but tedious.  It would be helpful to introduce diagrammatics for graded monoidal categories, similar to what was done for super 2-categories in \cite{EllisLauda_odd}.  We leave the details to the reader.
\end{proof}

\begin{example}\label{ex:koszulSign}
Let $\AC$ be a $\k$-linear monoidal category and $\AS=\KC^-(\AC)$ the homotopy category of bounded above complexes. As in Example \ref{ex:sign}, $\KC^-(\AC)$ has the the structure of a homologically graded monoidal category with grading shift defined by$A[1]_k=A_{k+1}$.  For any $A,B\in \KC(\AC)$, one can form the complex of homogeneous linear maps $\Homb_\AC(A,B)$, whose $k$-th homology group is isomorphic to the space of chain maps $A\rightarrow B[k]$ modulo homotopy.  In other words, $H(\Homb_\AC(A,B))\cong \Homg_{\KC(\AC)}(A,B)$.  The tensor product of morphisms defines a chain map (of complexes of $\k$-modules)
\[
\Homb_\AC(A,A')\otimes_\k \Homb_{\AC}(B,B')\rightarrow \Homb_\AC(A\otimes B, B\otimes B')
\]
Because of the signs involved in the tensor product of complexes, there must be signs involved in order that the above be a chain map.  The convention is to use the Koszul sign rule:
\[
(f\otimes g)_{A_i\otimes B_j} = (-1)^{i|g|}f_i\otimes g_j.
\]
It is an exercise to deduce from this rule the sign in (\ref{eq:signInTensorComp}).
\end{example}

If $E=\bigoplus_{i\in \Z} E_i$ is a graded ring, let $E^{\op}$ denote $E$ with the opposite multiplication $(\a,\b)\mapsto (-1)^{|\a||\b|}\b\cdot \a$.  A graded left $E$-module is a graded abelian group $M=\bigoplus_i M_i$ with a left $E$-action such that $\a M_i\subset M_{i+|\a|}$ for every homogeneous $\a\in E$.  A graded right $E$-module is the same thing as a graded left $E^{\op}$-module.  If $E$ is graded commutative, then any graded right $E$-module gives rise to a graded left $E$-module via $\a\cdot m = (-1)^{|\a||m|}m\cdot \a$ for all $\a\in E$, $m\in M$.

In the next proposition, let $\AS$ be a graded monoidal category.  For each $A\in \AS$, and each $\a\in \Endg_\AS(\one)$ of degree $k$, let $L_\a,R_\a\in \Endg_\AS(A)$ denote the compositions
\[
A \cong \one\otimes A\buildrel \a\otimes \Id_A\over \longrightarrow \one[k]\otimes A \cong A\otimes \one[k] \cong A[k]
\]
\[
A \cong A\otimes \one \buildrel \Id_A\otimes \a\over \longrightarrow A\otimes \one[k]\cong A[k]
\]
respectively.  
\begin{proposition}\label{prop:superActions}
The ring $\Endg\k\modules(\one)$ is graded commutative, and the maps $L,R:\Endg_\AS(\one)\rightarrow \Endg_\AS(A)$ satisfy
\begin{enumerate}
\item $L_\a\circ L_\b = L_{\a\circ \b}$.
\item $R_\a\circ R_\b = (-1)^{|\a||\b|}R_{\b\circ \a}$.
\item $L_\a\circ R_\b = (-1)^{|\a||\b|}R_{\b}\circ L_{\a}$.
\end{enumerate}
Furthermore, the four evident actions of $\Endg_\AS(\one)$ on $\Homg_\AS(A,B)$ are related by:
\begin{enumerate}
\item[(4)] $L_\a\circ f = (-1)^{|\a||f|} f\circ L_\a$.
\item[(5)] $R_\a\circ f = (-1)^{|\a||f|} f\circ R_\a$.
\end{enumerate}
for all $f\in \Homg_\AS(A,B)$ and all $\a\in \Endg_\AS(\one)$.  The maps $L_\a\circ (-)$ and $R_\a\circ (-)$ are called the left and right actions of $\Endg_\AS(\one)$ on $\Homg_\AS(A,B)$.  Finally, 
\begin{enumerate}
\item[(6)]
The composition of morphisms gives a map $\Homg_\AS(B,C)\otimes \Homg_\AS(A,B)\rightarrow \Homg_\AS(A,C)$ which is compatible with the left and right actions of $\Endg_\AS(\one)$.
\end{enumerate}
\end{proposition}
\begin{proof}
Straightforward.
\end{proof}

\subsection{Triangulated monoidal categories}
\label{subsec:triangMonoidal}
 We will be interested in categories $\AS$ which have the structures of a monoidal category and a triangulated category, and such that these structures are compatible in a certain way.  We adopt the axioms in \cite{MayTraces} for triangulated categories.  It seems that there is no general agreement on a set of axioms for triangulated monoidal categories, though in the symmetric monoidal case there are standard references \cite{MayTraces,Balmer10}.  For us, Definition \ref{def:triangMonoidal} is the most natural.

\begin{definition}\label{def:triang}
A \emph{triangulated category} is a category $\AS$ together with an invertible endofunctor $[1]:\AS\rightarrow \AS$, called the \emph{suspension}, and a collection $\TS$ of triangles
\[
A\buildrel f\over \rightarrow B\buildrel g\over \rightarrow C\buildrel h\over \rightarrow A[1],
\]
called \emph{distinguished triangles}, satisfying the following axioms:

\vskip7pt
\noindent \textbf{Axiom T1.} 
\begin{enumerate}
\item[(a)] Any triangle $A\buildrel \Id \over\rightarrow  A\rightarrow 0\rightarrow A[1]$ is distinguished.
\item[(b)] Any morphism $f:A\rightarrow B$ fits into a distinguished triangle $A\buildrel f\over\rightarrow B\rightarrow C\rightarrow A[1]$.
\item[(c)] Any triangle which is isomorphic to a distinguished triangle is distinguished.
\end{enumerate}

\vskip7pt
\noindent \textbf{Axiom T2.}  If $(f,g,h)$ is distinguished, then so is $(g,h,-f[1])$.  We say that $(g,h,-f[1])$ is obtained from $(f,g,h)$ by \emph{rotating}.

\vskip7pt
\noindent
\textbf{Axiom T3.}  Suppose we are given distinguished triangles
\[
\begin{diagram}A &\rTo^{f} & B & \rTo^{f'} & D & \rTo^{f''} & A[1]\end{diagram}
\]
\[
\begin{diagram}B &\rTo^{g} & C & \rTo^{g'} & E & \rTo^{g''} & B[1]\end{diagram}.
\]
Then for any choice of distinguished triangle
\[
\begin{diagram}
A &\rTo^{g\circ f} & C & \rTo^{h'} & F & \rTo^{h''} & A[1]\end{diagram}
\]
There exists a distinguished triangle
\[
\begin{diagram}D &\rTo^{j} & F & \rTo^{j'} & E& \rTo^{f'[1]\circ g''} & D[1]\end{diagram}
\]
such that the following diagram commutes:
\begin{equation}\label{eq:TateDiagram}
\begin{minipage}{5.1in}
\begin{tikzpicture}
\tikzstyle{every node}=[font=\small]
\node (a) at (0,0) {$A$};
\node (b) at (4,0) {$C$};
\node (c) at (8,0) {$E$};
\node (d) at (12,0) {$D[1]$};
\node (e) at (2,-2) {$B$};
\node (f) at (6,-2) {$F$};
\node (g) at (10,-2) {$B[1]$};
\node (h) at (4,-4) {$D$};
\node (i) at (8,-4) {$A[1]$};\\
\path[->,>=stealth',shorten >=1pt,auto,node distance=1.8cm, thick]
(a) edge node {$f$} (e)
(e) edge node {$f'$} (h)
(h) edge [bend right, looseness =1,out=-45,in=225] node[below] {$f''$} (i)
(i) edge node[label={[shift={(.7,-.7)}]$f[1]$}] {}  (g)
(g) edge node[label={[shift={(.7,-.7)}]$f'[1]$}] {}  (d)
(e) edge node[below right] {$g$} (b)
(b) edge [bend left, looseness=1,out=45,in=135]  node {$g'$} (c)
(c) edge node {$g''$} (g)
(c) edge [bend left, looseness=1,out=45,in=135] node {$f'[1]\circ g''$} (d)
(b) edge node {$h'$} (f)
(f) edge node {$h''$} (i)
(a) edge [bend left, looseness=1,out=45,in=135]  node {$g\circ f$} (b)
;
\path[dashed,->,>=stealth',shorten >=1pt,auto,node distance=1.8cm, thick]
(h) edge node[below right] {$j$} (f)
(f) edge node[below right] {$j'$} (c);
\end{tikzpicture}
\end{minipage}
\end{equation}
\end{definition}

It is a consequence of the axioms that the object $C$ in a distinguished triangle $A\buildrel f\over\to B\to C\to A[1]$ is uniquely determined by $f$ up to isomorphism.  We call $C$ the mapping cone on $f$, and we write $C\cong \Cone(f)$.  We have $\Cone(f)\cong 0$ iff $f$ is an isomorphism.  The axiom (T3) is called Verdier's axiom, and is a strengthened version of the octahedral axiom.  It states essentially that the mapping cone of a composition $A\buildrel f\over \rightarrow B\buildrel g\over \rightarrow C$ is isomorphic to the mapping cone on a map $h:\Cone(g)[-1]\rightarrow \Cone(f)$ in a particularly nice way.  Here $h$ is the composition of canonical maps
\[
\Cone(g)[-1]\rightarrow B\rightarrow \Cone(f).
\]

\begin{definition}\label{def:triangMonoidal}
A \emph{triangulated monoidal category} is a tuple $(\AS,\otimes,\one,\mu,\l,\rho,[1],\a,\b, \TS)$ such that $(\AS,[1],\TS)$ is a triangulated category and $(\AS,\otimes,\one,\mu,\l,\rho,[1],\a,\b)$ is a graded monoidal category of homological type  (see Definition \ref{def:gradedMonoidalCat}).   We further require that $\otimes$ be exact: that is, for any distinguished triangle $A\buildrel f\over \rightarrow B\buildrel g\over \rightarrow C\buildrel h\over\rightarrow A[1]$ and any $D\in \AS$ the following triangles are distinguished:
\begin{equation}\label{eq:DtimesT}
\begin{diagram}
D\otimes A & \rTo^{\Id_D\otimes f} & D\otimes B & \rTo^{\Id_D\otimes g} & D\otimes C & \rTo^{\phi_{D,A}\circ (\Id_D\otimes h)} & (D\otimes A)[1]
\end{diagram}
\end{equation}
\begin{equation}\label{eq:TtimesD}
\begin{diagram}
A\otimes D & \rTo^{f\otimes \Id_D} & B\otimes D & \rTo^{g\otimes \Id_D} & C\otimes D & \rTo^{\psi_{A,D}\circ (\Id_D\otimes h)} & (D\otimes A)[1]
\end{diagram}
\end{equation}
where $\phi_{D,A}:D\otimes A[1] \buildrel \cong\over\rightarrow (D\otimes A)[1]$ and $\psi_{A,D}:A[1]\otimes D \buildrel\cong\over\rightarrow (A\otimes D)[1]$ are the evident isomorphisms constructed from the natural transformations $\mu,\a,\b$.

This completes our definition of a triangulated monoidal category.  As usual, we will abuse language and simply refer to $\AS$ as a triangulated monoidal category.
\end{definition}

\begin{notation}\label{notation:triangleTimesObject}
We will denote the triangles (\ref{eq:DtimesT}) and (\ref{eq:TtimesD}) by $(\Id_D\otimes f, \Id_D\otimes g, \Id_D\otimes h)$ and $(f\otimes \Id_D, g\otimes \Id_D, h\otimes \Id_D)$, respectively, by abuse of notation.
\end{notation}
\begin{remark}
The fact that the sign $\pi$ in axiom (GM2) is required to be $-1$ for triangulated monoidal categories can be motivated informally in the following way.  Any triangle $(f,g,h)$ is isomorphic to $(-f,-g,h)$.  Rotating triangles three times implies that $(f[1],g[1],-h[1])$ is distinguished if $(f,g,h)$ is.  We would like the shift absorption isomorphism and the braiding morphism to yield an isomorphisms of triangles
\[
(\Id_{\one[1]}\otimes f, \Id_{\one[1]}\otimes g, \Id_{\one[1]}\otimes h) \ \ \ \ \cong \ \ \ \ \ (f[1],g[1],-h[1]) \ \ \ \ \cong \ \ \ \ (f\otimes\Id_{\one[1]},g\otimes \Id_{\one[1]},h\otimes\Id_{\one[1]}).
\]
The relative signs on the third maps forces $\pi=-1$.
\end{remark}

\begin{example}
If $\AC$ is an additive monoidal category, then $\KC^\pm(\AC)$ and $\KC^b(\AC)$ are triangulated monoidal categories in a canonical way.
\end{example}


\section{Complementary idempotents in monoidal categories}
\label{sec:firstProps}



\subsection{Main definitions}\label{subsec:definitions}
Throughout the remainder of this paper, let $\AS$ denote a triangulated monoidal category.

\begin{definition}\label{def-fixedCategories}
A \emph{weak idempotent} in $\AS$ is an object $\idempotent\in \AS$ such that $\idempotent^{\otimes 2}\cong \idempotent$.  If $\idempotent\in \AS$ is a weak idempotent, then let $\idempotent \AS\subset \AS$ denote the full subcategory consisting of objects $A\in \AS$ such that $A\cong\idempotent \otimes A$.  Similarly define $\AS \idempotent$ and $\idempotent \AS \idempotent:=\idempotent\AS \cap \AS\idempotent$. 
\end{definition}

Unfortunately, the subcategories $\idempotent \AS,\AS\idempotent,\idempotent\AS\idempotent\subset \AS$ are not triangulated in general (see Example \ref{ex:notCatIdempt}).  This suggests that the notion of weak idempotent is too general to be of any use.  A much better notion would involve morphisms in some way.  We will work exclusively with the following:

\begin{definition}\label{def:idempotents}
A \emph{pair of complementary idempotents} is a pair of objects $\counital,\unital\in \AS$ with $\unital\otimes \counital\cong \counital\otimes \unital\cong 0$, together with a distinguished triangle:
\begin{equation}\label{eq:PQtriangle}
\counital\buildrel \e\over \rightarrow \one \buildrel \eta\over\rightarrow \unital\buildrel \d\over \rightarrow \counital[1].
\end{equation}
We sometimes call any such distinguished triangle an \emph{idempotent triangle}.  We say $(\unital,\eta)$ is a \emph{unital idempotent} if it fits into an idempotent triangle; similarly, we say $(\counital,\e)$ is a \emph{counital idempotent}.  We will also write $\counital^c=\unital$ and $\unital^c = \counital$ for the complementary idempotents.
\end{definition}
In the presence of a triangulated monoidal structure, this definition of unital and counital idempotent agrees with the one provided in the first paragraph of the introduction of this paper.  In \S \ref{subsec-uniqunessOfIdempts} we will show that the distinguished triangle \ref{eq:PQtriangle} is determined by the isomorphism class $[\unital]$ (or $[\counital]$) up to unique isomorphism of triangles.  Thus for most purposes, the structure maps $\e,\eta,\d$ need not be explicitly included in the definition.  

\begin{observation} If $(\unital,\eta)$ is a unital idempotent, then $\unital^c\cong \Cone(\eta)[-1]$.  If $(\counital,\e)$ is a counital idempotent, then $\counital^c\cong \Cone(\e)$.
 \end{observation}

\begin{theorem}\label{thm-pApMonoidal}
Let $\idemp$ be a unital or counital idempotent in $\AS$.  Then $\idemp\AS \idemp$ has the structure of a monoidal category with monoidal identity $\idemp$.\qed
\end{theorem}
As was pointed out in Remark 2.19 of \cite{BD14} this is false for weak idempotents.  
\begin{proof}
This is Lemma {2.18} in \cite{BD14}.
\end{proof}

We now generalize the notion of unital and counital idempotent:

\begin{definition}\label{def:catIdempts}
Let $\AS$ be a triangulated monoidal category, and let $\idempotent\in \AS$ be a weak idempotent.  A \emph{unital idempotent relative to $\idempotent$} is an object $\unital\in \AS$ together with a map $\eta: \idempotent\rightarrow \unital$ such that $\eta\otimes\Id_\unital:\idempotent \otimes \unital\rightarrow \unital\otimes \unital$ and $\Id_\unital\otimes \eta:\unital\otimes \idempotent \rightarrow \unital\otimes \unital$ are isomorphisms.  One has a similar notion of \emph{counital idempotent relative to $ \idempotent$}.

We say $\idemp\in\AS$ is a \emph{categorical idempotent} if there is a chain of objects $\one=\idempotent_0,\ldots, \idempotent_r=\idempotent$ such that $\idempotent_i$ has the structure of a unital or counital idempotent relative to $\idempotent_{i-1}$ for $1\leq i\leq r$.    The \emph{depth} of a categorical idempotent is the minimal $r$ in such a chain.
\end{definition}

\begin{remark}\label{rmk:depth}
The only catetgorical idempotent of depth 0 is $\one$.  A categorical idempotent has depth $1$ if and only if it has the structure of a unital or counital idempotent.
\end{remark}

\begin{remark}\label{rmk:depth 2}
It is possible to show that if $\unital$ and $\counital$ are a unital and a counital idempotent which commute, then $\idempotent:=\unital\otimes \counital$ is a categorical idempotent of depth 2 (or less) .   Proposition \ref{prop:classificationOfCatIdempts} implies that every categorical idempotent arises in this way, so that every categorical idempotent has depth $\leq 2$.
\end{remark}
\begin{remark}\label{rmk:relIdemptsMonoidal}
Theorem \ref{thm-pApMonoidal} implies that if $\idemp$ is a categorical idempotent, then $\idempotent \AS\idempotent$ has the structure of a monoidal category with monoidal identity $\idempotent$.
\end{remark}
\begin{remark}\label{rmk:relIdempts}
An object $\otherIdempotent\in \AS$ has the structure of a (co)unital idempotent relative to $\idempotent$ if and only if $\otherIdempotent$ is a (co)unital idempotent in $\idempotent \AS\idempotent$
\end{remark}

Our aim now is to show that if $\idempotent$ is a categorical idempotent, then  $\idempotent \AS, \AS\idempotent,\idempotent\AS\idempotent\subset \AS$ are full triangulated subcategories.  As a first step we have the following:
\begin{proposition}[Projector absorbing]\label{prop-projectorAbsorbing}
Let $\AS$ be a triangulated monoidal category, and $(\counital,\unital)$ a pair of complementary idempotents in $\AS$.  For $A\in \AS$, The following are equivalent:
\begin{enumerate}
\item $\unital\otimes A\cong A$.
\item $\eta\otimes\Id_A:\one\otimes A\rightarrow \unital\otimes A$ is an isomorphism.
\item $\counital\otimes A\cong 0$.
\end{enumerate}
A similar statement holds with ``left'' replaced by ``right,'' or with the roles of $\unital$ and $\counital$ reversed.
\end{proposition}
Note that this proposition allows us to transfer information about objects (e.g.~$\counital\otimes A\cong 0$) to information about morphisms (e.g.~$\eta\otimes \Id_A:\one\otimes A\rightarrow \unital\otimes A$ is an isomorphism).
\begin{proof}
We consider only the case when $A$ absorbs the unital idempotent $\unital$ from the left.  The other cases are similar.  Clearly (2) implies (1).  Suppose now that $\unital\otimes A\cong A$, then
\[
\counital\otimes A \cong \counital\otimes (\unital\otimes A)\cong (\unital\otimes \counital)\otimes A\cong 0
\]
This shows that (1) implies (3).

Now, assume $\counital\otimes A\cong 0$.  Then tensoring the triangle $\counital\rightarrow \one\rightarrow \unital\rightarrow \counital[1]$ on the right with $A$ gives a distinguished triangle
\[
0\rightarrow \one\otimes A\buildrel \eta\otimes \Id_A\over \longrightarrow \unital\otimes A\rightarrow 0
\]
A standard fact about triangulated categories implies that $\eta\otimes \Id_A$ is an isomorphism.  This completes the shows that (3) implies (2), and completes the proof.
\end{proof}

\begin{proposition}\label{prop-coneClosure}
If $\idemp\in\AS$ is categorical idempotent, then $\idemp\AS,\: \AS \idemp,\: \idemp\AS \idemp\subset \AS$ are closed under mapping cones and the shift $[1]$.
\end{proposition}
\begin{proof}
Recall the notion of depth of a categorical idempotent (Definition \ref{def:catIdempts}). Any easy induction on depth reduces to the case when $\idempotent$ has depth 1 (i.e.~is a unital or counital idempotent).  So assume $\idempotent$ has depth 1.  Let $\idempotent^c$ be the complementary idempotent.  Suppose we have a distinguished triangle
\begin{equation}\label{eq-aTriangle1}
A\rightarrow B\rightarrow C\rightarrow A[1]
\end{equation}
in $\AS $.  If $A,B\in \idemp\AS $, then $\idempotent^c\otimes A\cong \idempotent^c\otimes B\cong 0$  by Proposition \ref{prop-projectorAbsorbing}.  Then tensoring the triangle (\ref{eq-aTriangle1}) on the left with $\otherIdemp$ gives a distinguished triangle
\[
0\rightarrow 0\rightarrow \idempotent^c\otimes C\rightarrow 0
\]
Thus $\idempotent^c\otimes C\cong 0$.  By Proposition \ref{prop-projectorAbsorbing} again, this implies that $C\in \idemp\AS$.  That is to say, if two terms of a distinguished triangle are in $\idemp\AS $, then so is the third.  A similar argument takes care of the cases of $\AS \idemp$ and $\idemp\AS \idemp$.\end{proof}

\begin{example}\label{ex:notCatIdempt}
Let $(\counital,\unital)$ be a pair of complementary idempotents in $\AS$ such that $\one\not\cong \unital\oplus \counital$.  Then $\idempotent :=\unital\oplus \counital$ is not a categorical idempotent, though it is a weak idempotent: $\idempotent\otimes \idempotent \cong \idempotent$.  One can see this by showing that $\idempotent \AS$ is not triangulated.  Indeed $\idempotent \otimes \unital\cong \unital$ and $\idempotent\otimes \counital\cong \counital$, but $\idempotent\otimes \one$ is not isomorphic to $\one$, by hypothesis.  On other hand there is a distinguished triangle:
\[
\counital\rightarrow \one\rightarrow \unital \rightarrow \counital[1]
\]
The first and third terms are in $\idempotent \AS$, but the middle term is not, hence $\idempotent \AS$ is not triangulated.
\end{example}

\subsection{Semi-orthogonality properties}
\label{subsec-semiOrtho}

In this section we show that any pair of complementary idempotents enjoys a certain orthogonality relation with respect to $\Hom_{\AS }(-,-)$.

\begin{theorem}\label{thm-semiOrtho}
Let $(\counital,\unital)$ be a pair of complementary idempotents in $\AS $.  Let $A,B\in \AS$ be arbitrary.  Then
\[
\Hom_{\AS}(\counital\otimes A,\unital\otimes B) \cong 0 \cong \Hom_{\AS}(A\otimes \counital,B\otimes \unital).
\]
\end{theorem}
As a special case, we see that $\Homg_{\AS }(\counital,\unital)\cong 0$.  Note that homs in the other direction may be nonzero.  For instance, by definition there is a map $\d:\unital\rightarrow \counital[1]$ such that $\one\cong \Cone(\d)[-1]$.  If $\d=0$, then this would mean that $\one\cong \unital\oplus \counital$, which is not generally true.
\begin{proof}
Assume that $\counital\otimes A\cong A$ and $\unital\otimes B\cong B$.  We will prove that $\Hom_{\AS}(A,B)\cong 0$.  Since $\one$ is the monoidal identity, we have an isomorphism
\[
\Phi:\Hom_{\AS }(A,B)\rightarrow \Hom_{\AS }(\one\otimes A, \one\otimes B)
\]
which sends $f\mapsto \Id_\one\otimes f$.  By hypothesis $ \counital\otimes A \cong A$ and $\unital\otimes B \cong B$.  Then Proposition \ref{prop-projectorAbsorbing} says that $\e\otimes \Id_{A}:\counital \otimes A \rightarrow \one \otimes A$ and $\eta\otimes \Id_{B}:\one \otimes B \rightarrow \unital\otimes B$ are isomorphisms in $\AS $.  Thus we have an isomorphism
\[
\Psi:\Hom_{\AS }(\one\otimes A, \one\otimes B)\rightarrow \Hom_{\AS }(\counital\otimes A, \unital\otimes B)
\]
which sends $g\mapsto (\eta\otimes \Id_{B})\circ g\circ (\e\otimes \Id_{A})$.  The composition $\Psi\circ \Phi$ sends $f\in \Hom_{\AS }(A,B)$ to
\[
(\eta\otimes \Id_{B})\circ (\Id_\one\otimes f)\circ (\e\otimes\Id_{A}) = (\Id_\unital\otimes f)\circ ((\eta\circ \e)\otimes \Id_{A})
\]
This latter map is zero since $\eta$ and $\e$ are adjacent maps in a distinguished triangle.  So on one hand $\Psi\circ \Phi$ is an isomorphism $\Hom_{\AS }(A,B)\cong \Hom_{\AS}(\counital \otimes A,\unital \otimes B)$, but on the other hand $\Psi\circ \Phi$ is the zero map.  Thus $\Hom_{\AS }(A,B)\cong 0$.  A similar argument takes care of the case when $A\otimes \counital \cong A$ and $B\otimes \unital\cong B$.
\end{proof}

\begin{remark}\label{rmk-semiOrthoCx}
If $\AS $ is the homotopy category of complexes over a monoidal category $\CS$, then $\Hom_{\CS}(A,B)$ is the zeroth homology of the hom complex $\Homb_{\CS}(A,B)$, for all complexes $A,B$.  If $A\in \AS  \counital$ and $B\in \AS  \unital$, then the above proof can be modified into a proof that $\Homb_{\CS}(A,B)\simeq 0$.
\end{remark}

\begin{remark}\label{rmk:semi-ortho}
Let $\MS$ be a triangulated category on which $\AS$ acts by exact endofunctors, and let $(\counital,\unital)$ be a pair of complementary idempotents.  Then for any $M\in \MS$, there is a distinguished triangle $\counital(M)\rightarrow M \rightarrow \unital(M)\rightarrow \counital(M)[1]$, and the above proof can be modified in an obvious way to prove that $\Homg_\MS(\counital(M),\unital(M))\cong 0$ for all $M$.
\end{remark}

The following is a very useful consequence of the semi-orthogonality relation between an idempotent and its complement.
\begin{proposition}\label{prop-resProp}
Suppose $(\unital,\eta)$ is a unital idempotent in a triangulated monoidal category $\AS $.  Let $B\in\AS \unital$ and $A\in \AS $ be arbitrary.  Then precomposition with $\Id_A\otimes \eta$ gives an isomorphism
\[
\Hom_{\AS }(A\otimes \unital, B)\cong \Hom_{\AS }(A\otimes \one, B).
\]
Similarly, if $(\counital,\e)$ is a counital idempotent and $A\in \AS \counital$, $B\in\AS $, then post-composition with $\Id_B\otimes \e$ defines an isomorphism
\[
\Hom_{\AS }(A,B\otimes \counital)\cong \Hom_{\AS }(A,B\otimes \one).
\]
\end{proposition}
\begin{proof}
Let us prove the first statement.  The second is similar.  Suppose $B\in\AS \unital$, and let $A\in\AS $ be arbitrary.  By Theorem \ref{thm-semiOrtho} we have $\Hom_{\AS }(A\otimes \counital, B)\cong 0$.  Consider the distinguished triangle
\[
\counital\rightarrow \one\rightarrow \unital\rightarrow \counital[1]
\]
Tensoring with $A$ gives a distinguished triangle
\[
A\otimes \counital\rightarrow A\otimes \one \rightarrow A\otimes \unital\rightarrow A\otimes \counital[1]
\]
Since $\AS $ is a triangulated category, applying the functor $\Hom_{\AS }(-,B)$ to the above distinguished triangle yields a long exact sequence, in which the relevant part is
\[
\Hom_{\AS }(A\otimes \counital[1], B)\rightarrow \Hom_{\AS }(A\otimes \unital, B)\rightarrow \Hom_{\AS }(A\otimes \one,B)\rightarrow \Hom_{\AS }(A\otimes \counital,B)
\]
The two extremal terms are zero by Theorem \ref{thm-semiOrtho}, hence the middle map is an isomorphism.  This map is precisely precomposition with $\Id_A\otimes \eta$.  This completes the proof.
\end{proof}

\begin{remark}\label{rmk:adjoint}
Let $\AS$ act on $\MS$, and let $(\counital,\unital)$ be a pair of complementary idempotents.  The proof of Proposition \ref{prop-resProp} can be modified in an obvious way to give a natural isomorphism $\Hom_\MS(\unital(M),N)\cong \Hom_\MS(M,N)$ for every $N\in \im \unital$.  This implies that the inclusion $\im \unital\rightarrow \MS$ has a left adjoint given by the functor $M\rightarrow \unital(M)$.  Conversely, any fully faithful functor $\NS\rightarrow \MS$ with a right adjoint determines a unital idempotent endofunctor of $\MS$ (an exercise in unpacking definitions).   Similarly, $N\mapsto \counital(N)$ determines a right adjoint to the inclusion $\im \counital \rightarrow \MS$, and any fully faithful functor with a right adjoint determines a counital idempotent endofunctor.
\end{remark}
\begin{corollary}\label{cor-endRings}
Precomposition with $\one\buildrel\eta\over\rightarrow \unital$ gives an isomorphism $\Homg(\unital,\unital)\rightarrow \Homg(\one,\unital)$.  Post-composition with $\counital\buildrel\e\over\rightarrow \one$ gives an isomorphism $\Homg(\counital,\counital)\rightarrow \Homg(\counital,\one)$.\qed
\end{corollary}


We now show that unital idempotents are also unital algebras.  Similar arguments show that counital idempotents are counital coalgebras.

\begin{lemma}\label{lemma-unitsAreSame}
We have $\eta\otimes \Id_\unital = \Id_\unital\otimes \eta$, regarded as maps $\unital\rightarrow \unital\otimes \unital$.  
\end{lemma}
\begin{proof}
By composing with the unitors in $\AS $, we will regard $\eta\otimes \Id_\unital$ and $\Id_\unital\otimes \eta$ as elements of $\Hom(\unital,\unital\otimes \unital)$.  By Proposition \ref{prop-resProp}, precomposing with $\eta$ gives an isomorphism $\eta^{\ast}:\Hom(\unital,\unital\otimes \unital)\rightarrow \Hom(\one,\unital\otimes \unital)$.  Clearly
\[
\eta^{\ast}(\eta\otimes \Id_\unital)=\eta\otimes \eta=\eta^{\ast}(\Id_\unital\otimes \eta)
\]
Since $\eta^{\ast}$ is an isomorphism, we conclude that $\eta\otimes \Id_\unital=\Id_\unital\otimes \eta$.
\end{proof}
Define $\mu:\unital\otimes \unital\buildrel\cong\over\rightarrow \unital$ to be an inverse to $\eta\otimes \Id_\unital = \Id_\unital\otimes \eta$.
\begin{proposition}\label{prop-pAsAlgebra}
The isomorphism  $\mu:\unital\otimes \unital\buildrel \cong \over \rightarrow \unital$ makes $\unital$ an associative algebra in $\AS $ with unit $\eta:\one\rightarrow \unital$.
\end{proposition}
\begin{proof}
The unit axiom holds since by definition $\mu\circ (\eta\otimes \Id_\unital) = \Id_\unital$ and $\mu\circ(\Id_\unital\otimes \mu)=\Id_\unital$ (after the insertion of the left and right unitors in $\AS $, where appropriate).  For associativity, Proposition \ref{prop-resProp} says that precomposition with $\eta^{\otimes 3}$ gives an isomorphism $(\eta^{\otimes 3})^\ast:\Hom(\unital^{\otimes 3},\unital)\rightarrow \Hom(\one^{\otimes 3},\unital)$.  The unit axiom implies that 
\[
\mu\circ (\Id_\unital\otimes \mu)\circ \eta^{\otimes 3} = \mu\circ(\mu\otimes \Id_\unital)\circ \eta^{\otimes 3} = \eta
\]
from which associativity follows.
\end{proof}


\subsection{Endomorphism algebras}
\label{subsec-endIdempts}
Recall that a graded ring $R=\bigoplus_{i\in \Z} R$ is graded commutative if $ab=(-1)^{|a||b|}ba$ for all $a,b\in R$.  The graded center of a graded ring is the ring $Z(R)\subset R$ spanned of homogeneous elements $z\in R$ such that $za=(-1)^{|z||a|}$ for all $a\in R$.  If $\k$ is a graded commutative ring, then a graded $\k$-algebra is a graded ring $R$ with a homogeneous map $\k\rightarrow Z(R)$.  In particular scalars \emph{super-commute} with elements of $R$.

The following is easy but useful.  

\begin{theorem}\label{thm-endRings}
Let $\idempotent$ be a categorical idempotent in a triangulated monoidal category $\AS$.  Then $\Endg_\AS(\idempotent)$ is a graded commutative $\Endg_\AS(\one)$-algebra.  The unit map $\Endg_\AS(\one)\rightarrow \Endg_\AS(\idemp)$ can be defined in one of two ways, which coincide:
\[
\a \mapsto \a\otimes \Id_{\idempotent} = \Id_{\idempotent}\otimes \a.
\]
That is, $R_\a = L_\a :\Endg_\AS(\one)\rightarrow \Endg_\AS(\idempotent)$, where $L$ and $R$ are defined in the remarks preceding Proposition \ref{prop:superActions}.
\end{theorem}
\begin{proof}
We can regard $\idempotent$ as the monoidal identity in a graded monoidal category $(\idempotent \AS  \idempotent, \otimes ,\idempotent)$, so $\Endg_\AS(\idempotent)$ is graded commutative by Proposition \ref{prop:superActions}.

We will now prove the second statement.  An easy induction reduces the general case to the case where $\idempotent$ is a unital or counital idempotent. Assume first that $\idempotent$ is a unital idempotent with unit map $\eta:\one\rightarrow \idemp$.  Then precomposition with $\eta$ gives an isomorphism
\[
\Endg_\AS(\idemp)\rightarrow \Homg_\AS(\one,\idemp).
\]
This map is compatible with the left and right actions of $\Endg_\AS(\one)$.  On the other hand, these actions coincide up to sign for $\Homg_\AS(\one,\idemp)$, since they can be described in terms of the left and right actions on the first argument $\one$.  Thus, the left and right actions of $\Endg_\AS(\one)$ on $\Endg_\AS(\idemp)$ coincide up to sign.  In particular $\a\otimes \Id_\idemp = \Id_{\idemp}\otimes \a$; there is no sign since $\Id_\idemp$ has degree zero.  A similar argument takes care of the case when $\idempotent$ is a counital idempotent.  This completes the proof.  
\end{proof}

\subsection{Commuting properties}
\label{subsec:commutingProps}
In this section we show that an object $A\in \AS$ commutes with $\unital$ if and only if it commutes with the complement $\unital^c$.  Further, the full category consisting of such objects $A$ is triangulated:

\begin{proposition}\label{prop:commuteWithidemp}
Let $(\counital, \unital)$ be a pair of complementary idempotents in $\AS$, and let $A\in \AS$ be arbitrary.  The following are equivalent:
\begin{enumerate}
\item $\unital\otimes A\cong A\otimes \unital$.
\item $\counital\otimes A\cong A\otimes \counital$.
\item $\counital\otimes A\otimes \unital \cong 0 \cong \unital\otimes A\otimes \counital$.
\end{enumerate}
Furthermore, if $\BS\subset \AS$ denotes the full subcategory consisting of objects $A$ which satisfy either one of these equivalent conditions, then $\BS$ is a triangulated subcategory. 
\end{proposition}
\begin{proof}
Suppose (1) holds.  Then
\[
\counital\otimes A \otimes \unital \cong \counital\otimes \unital \otimes A \cong  0,
\]
since $\unital$ and $\counital$ are orthogonal.  A similar argument shows $\unital\otimes A\otimes \counital\cong 0$.  Thus, (1) implies (3).  A similar argument shows that (2) implies (3).

The implications (3) $\Rightarrow$ (1) and (3)$\Rightarrow$ (1) follow from projector absorbing (Proposition \ref{prop-projectorAbsorbing}).  This proves the equivalence of (1),(2),(3).  The last statement is obvious, since condition (3) is clearly closed under mapping cones and suspension $[1]$.  This completes the proof.
\end{proof}

\begin{theorem}\label{thm-centralIdempts}
Let $\unital \in \AS $ be a unital idempotent and $\counital$ its complement. The following are equivalent:
\begin{enumerate}
\item[(1)] $\unital\AS =\AS \unital$.
\item[(2)] $\unital\otimes A\cong A\otimes \unital$ for all $A\in \AS$.
\item[(3)] $\unital$ is central, i.e.~ there is a natural isomorphism $\unital\otimes A\cong \unital\otimes A$ for all $A\in\AS $.
\item[(4)] $\counital \AS  = \AS  \counital$.
\item[(5)] $\counital \otimes A\cong A\otimes \counital$ for all $A\in \AS$.
\item[(6)] $\unital$ is central, i.e.~ there is a natural isomorphism $\unital\otimes A\cong \unital\otimes A$ for all $A\in\AS $.
\end{enumerate}
\end{theorem}
\begin{proof}
Clearly (3)$\Rightarrow$(2) and (2)$\Rightarrow$ (1).  To show (1)$\Rightarrow$ (3), assume that $\unital\AS  = \AS  \unital$ (hence $\unital\AS = \unital\AS  \unital$), and let $A\in \AS $ be arbitrary.  Then $\unital\otimes A\in \unital\AS =\unital \AS  \unital$, and $A\otimes \unital\in \AS  \unital = \unital\AS  \unital$.  Since $\unital\AS  \unital$ has the structure of a monoidal category with monoidal identity $\unital$, we have isomorphisms
\[
 \unital\otimes A \rightarrow (\unital\otimes A)\otimes \unital \rightarrow \unital\otimes (A\otimes \unital) \rightarrow A\otimes \unital.
 \]
which are natural in $A$.  The first and last isomorphisms are the right and left unitors in $\unital \AS  \unital$.  This proves the equivalence of (1)-(3).  A similar argument proves the equivalence of (4)-(6).  The equivalence (2)$\Leftrightarrow$(5) is an immediate consequence of Proposition \ref{prop:commuteWithidemp}.
\end{proof}

\subsection{Fundamental theorem of categorical idempotents}
\label{subsec-uniqunessOfIdempts}

Consider the following example from linear algebra: let $A=M_{2\times 2}(\Z)$ be the algebra of 2-by-2 matrices, and consider
\[
p := \matrix{1&0\\0&0} \ \ \ \ \text{ and } \ \ \ \ p' := \matrix{1&1\\0&0}
\]
Then $p$ and $p'$ are idempotents, and satisfy $pp'=p'$, $p'p=p$.  In this section we see that the categorical analogue of  this left/right asymmetry never occurs.

\begin{theorem}\label{thm:idemptOrder}
Let $(\unital_i,\eta_i)$ be unital idempotents, which fit into distinguished triangles
\begin{equation}\label{eq-aPairOfResolutions}
\begin{diagram}[small]
 \counital_i & \rTo^{\e_i} & \one & \rTo^{\eta_i} & \unital_i& \rTo^{\d_i} &\counital_i[1]
\end{diagram}
\end{equation}
for $i=1,2$.  The following are equivalent:
\vskip7pt
\begin{minipage}{3in}
\begin{enumerate}
\item $\unital_1\otimes \unital_2\cong \unital_2$.
\item $\unital_2\otimes \unital_1\cong \unital_2$.
\item $\exists \nu:\unital_1\rightarrow \unital_2$ such that $\nu\circ \eta_1 = \eta_2$.
\end{enumerate}
\end{minipage}
\begin{minipage}{3in}
\begin{itemize}
\item[(4)] $\counital_1\otimes \counital_2\cong \counital_1$.
\item[(5)] $\counital_2\otimes \counital_1\cong \counital_1$.
\item[(6)] $\exists \theta:\counital_1\rightarrow \counital_2$ such that $\e_2\circ \theta = \e_1$.
\end{itemize}
\end{minipage}
\vskip7pt
If any one of these equivalent conditions is satisfied, then the maps $\nu,\theta$ from (3) and (6) are unique, so that there is a unique map of triangles:
\begin{equation}\label{eq-mapOfResolutions}
\begin{diagram}
 \counital_1 & \rTo^{\e_1} & \one & \rTo^{\eta_1} & \unital & \rTo^{\d_1} & \counital[1]\\
\dTo^{\theta}&& \dTo^{\Id} && \dTo^{\nu}&& \dTo^{\theta[1]}\\
 \counital_2 & \rTo^{\e_2} & \one & \rTo^{\eta_2} & \unital_2 & \rTo^{\d_2} & \counital_2[1]
\end{diagram}
\end{equation}
Furthermore $\nu$ realizes $\unital_2$ as a unital idempotent relative to $\unital_1$, while $\theta$ realizes $\counital_1$ as a counital idempotent relative to $\counital_2$.  
\end{theorem}
\begin{proof}
For the equivalence (1)$\Leftrightarrow$(4), two applications of projector absorbing (Proposition \ref{prop-projectorAbsorbing}) give
\[
\unital_1\otimes \unital_2\cong \unital_2 \ \ \ \ \Leftrightarrow  \ \ \ \ \counital_1\otimes \unital_2\cong 0 \ \ \ \ \Leftrightarrow  \ \ \ \   \counital_1\otimes \counital_2\cong \counital_1.
\]
A similar argument gives the equivalence (2)$\Leftrightarrow$(5).   We now prove that (1)$\Leftrightarrow$(3).

(3)$\Rightarrow$(1):  Suppose $\nu:\unital_1\rightarrow \unital_2$ is such that $\nu \circ \eta_1 = \eta_2$.  Let $\phi:\unital_1 \otimes \unital_2 \rightarrow \one \otimes \unital_2$ denote the composition
\[
\unital_1\otimes \unital_2 \buildrel \nu \otimes \Id_{\unital_2}\over \longrightarrow \unital_2\otimes \unital_2\buildrel (\eta_2\otimes \Id_{\unital_2})\inv\over\longrightarrow \one\otimes \unital_2
\]
Then set $\psi:=(\eta_1\otimes \Id_{\unital_2}):\one\otimes \unital_2\rightarrow \unital_1\otimes \unital_2$.  We claim that $\psi$ and $\phi$ are mutually inverse isomorphisms.  First, note that
\[
\phi\circ \psi = (\eta_2\otimes \Id_{\unital_2})\inv\circ ((\nu\circ \eta_1)\otimes \Id_{\unital_2}) = \Id_{\one \otimes \unital_2}
\]
since $\nu\circ \eta_1 = \eta_2$.  On the other hand, precomposition with $\psi$ gives an isomorphism
\[
\psi^\ast:\Hom_{\AS}(\unital_1\otimes \unital_2, \unital_1\otimes \unital_2)\rightarrow \Hom_{\AS }(\one\otimes \unital_2,\unital_1\otimes \unital_2)
\]
by Proposition \ref{prop-resProp}.  The equality $\phi\circ \psi = \Id_{\one\otimes \unital_2}$ ensures that $\Id_{\unital_1\otimes \unital_2}$ and $\psi\circ \phi$ have the same image under $\psi^\ast$, hence $\psi\circ \phi=\Id_{\unital_1\otimes \unital_2}$.  This proves that $\unital_1\otimes \unital_2\cong \unital_2$, and we conclude that (5)$\Rightarrow$(3).

(1)$\Rightarrow$(3):  Suppose that $\unital_1\otimes \unital_2\cong \unital_2$.  Then Proposition \ref{prop-resProp} says that precomposition with $\eta_1:\one\rightarrow \unital_1$ gives an isomorphism $\Hom_{\AS }(\unital_1,\unital_2)\cong \Hom_{\AS }(\one,\unital_2)$.  In particular there is a unique $\nu:\unital_1\rightarrow \unital_2$ such that $\nu\circ \eta_1=\eta_2$.  Thus, (1)$\Leftrightarrow$(3).  Similar arguments prove that (2)$\Leftrightarrow$(3), and (4)$\Leftrightarrow$(5)$\Leftrightarrow$(6).  This proves that statements (1)-(6) are equivalent, as desired.

We have seen above that $\nu$ is unique when it exists.  The proof that (1)$\Leftrightarrow$(3) above actually shows that $\nu$ becomes an isomorphism after tensoring with $\unital_2$ on the left (and by a symmetric argument, on the right).  Thus, $(\unital_2,\nu)$ is a unital idempotent in $\unital_1\AS  \unital_1$.  A similar argument establishes that $\theta:\counital_1\rightarrow \counital_2$ from (6) is unique, and that $(\counital_1,\theta)$ is a counital idempotent in $\counital_2\AS  \counital_2$. This completes the proof.
\end{proof}
In \cite{BD14} a similar result is obtained (with appropriate modifications) without the technology of triangulated categories. Their result is not much more general, since any (additive) monoidal category $\MC$ embeds monoidally in $\K^b(\MC)$.   As was observed in \emph{loc.~cit.}, the above proposition allows us to put a partial order on the set of unital (respectively counital) idempotents.

\begin{definition}\label{def:partialOrder}
Put a partial order on the set of (isomorphism classes of) unital idempotents in $\AS$ by declaring that $\unital_2\leq \unital_1$ if $\unital_1\otimes \unital_2\cong \unital_2$.
\end{definition}

\begin{remark}
Similar definitions apply to counital idempotents, so that the mapping $\unital\mapsto \unital^c$ defines an order-reversing bijection between between unital and counital idempotents.
\end{remark}

The following is often useful.
\begin{corollary}\label{cor:unitalCounitalOrtho}
Let $\unital\in \AS$ and $\counital\in \AS$ be a unital and counital idempotent.  Then $\unital\otimes \counital\cong 0$ if and only if  $\counital \otimes \unital \cong 0$.  If these are satisfied, then $\Homg_\AS(X,Y)\cong 0$ for any $X\in \counital\AS$ (respectively $X\in \AS\counital$) and any $Y\in \unital\AS$ (respectively $Y\in \AS\unital$).
\end{corollary}
\begin{proof}
Let $\counital^c$ denote the complementary idempotent to $\counital$.  Then
\[
\counital\otimes \unital \cong 0  \Leftrightarrow \counital^c\otimes \unital \cong \unital
\]
which holds if and only if $\unital \leq \counital^c$ (by definition).  Similar considerations imply that $\unital\otimes \counital$ iff $\unital\leq \counital^c$.  This proves the first statement.  The statements regarding the vanishing of $\Homg_{\AS}(X,Y)$ now follow from Theorem \ref{thm-semiOrtho}.
\end{proof}

As another corollary we have our uniqueness theorem.

\begin{theorem}\label{thm-uniquenessOfIdempts}
Let $(\counital_1,\unital_1)$ and $(\counital_2,\unital_2)$ be two pairs of complementary idempotents which fit into triangles (\ref{eq-aPairOfResolutions}).  The following are equivalent:
\begin{enumerate}
\item[(i)] $\unital_1\AS  = \unital_2\AS $
\item[(ii)] $\counital_1\AS  = \counital_2\AS  $
\item[(iii)] $\AS  \unital_1 = \AS  \unital_2$
\item[(iv)] $\AS  \counital_1 = \AS  \counital_2$
\end{enumerate}
Any one of the above equivalent statements implies that $\counital\cong \counital_2$ and $\unital\cong \unital_2$ canonically.  In fact, the canonical map of triangles (\ref{eq-mapOfResolutions}) is an isomorphism.
\end{theorem}
\begin{proof}
The equivalences (i)-(iv) follow from Theorem \ref{thm:idemptOrder}.  Uniqueness of the resulting map of triangles (\ref{eq-mapOfResolutions}) implies that this map is an isomorphism of triangles, as claimed.
\end{proof}
The following is a restatement of the uniqueness theorem which is often useful in practice.  In the statement, a \emph{left tensor ideal} is a full subcategory $\IC\subset \AS$ which is closed under $A\otimes -$ for all $A\in \AS$ (and analogously for right tensor ideal).

\begin{corollary}\label{cor-appliedUniqueness}
Let $\IC_L,\IC_R\subset \AS $ be a left and right tensor ideal, respectively, in a triangulated monoidal category $\AS$.  Let $\JC_R\subset \AS$ (respectively $\JC_L\subset \AS$) be the full subcategory of objects $J$ such that $I\otimes J\cong 0$ for all $I\in \IC_L$ (respectively $J\otimes I\cong 0$ for all $I\in \IC_R$).  Suppose there is an object $\unital\in \IC_L\cap \IC_R$ and a map $\eta:\one\rightarrow \unital$ such that $\Cone(\eta)\in \JC_L\cap \JC_R$.  Then
\begin{enumerate}
\item $(\unital,\eta)$ is a unital idempotent in $\AS $ such that $\unital\AS = \IC_R$ and $\AS\unital=\IC_L$.  In particular, the following are equivalent:
\begin{enumerate}
\item $\IC_L=\IC_R$.
\item $\JC_L=\JC_R$.
\item $\unital$ is central in $\AS$.
\end{enumerate}
\item $(\unital,\eta)$ is unique up to canonical isomorphism in the following way: if $(\unital',\eta')$ is another pair with these properties then there is a unique map $\phi:\unital\rightarrow \unital'$ such that $\phi\circ \eta = \eta'$, and this map is an isomorphism.
\end{enumerate}
Similar statements hold for counital idempotents, with the obvious modifications.
\end{corollary}
\begin{proof}
Let $\unital\in \IC_L\cap \IC_R$ be given, and suppose $\eta:\one\rightarrow \unital$ satisfies $\Cone(\eta)\in\JC_L\cap \JC_R$.  Then $\unital\otimes \Cone(\eta)\cong  0 \cong \Cone(\eta)\otimes \unital$, which implies that $(\unital,\eta)$ is a unital idempotent with complement $\unital^c:=\Cone(\eta)[-1]$. 

We claim that $A\in \AS $ satisfies $A\otimes \unital\cong A$ if and only if $A\in \IC_R$.  If $A\otimes \unital\cong A$, then $\AS \in\IC_L$ since $\unital\in \IC_L\cap \IC_R$ and $\IC_L$ is a left tensor ideal.  Conversely, if $A\in \IC_L$, then $A\otimes \unital^c\cong 0$ since $\unital^c \in \JC_L\cap \JC_R$, which implies that $\unital\otimes A\cong A$ by projector absorbing.   This proves (1).  Statement (2) follows immediately from the uniqueness theorem (Theorem \ref{thm-uniquenessOfIdempts}).  The statement about centrality follows immediately from Theorem \ref{thm-centralIdempts}.
\end{proof}

\section{Postnikov systems and decompositions of identity}
\label{sec:decompOfOne}
This section essentially concerns what might be called the categorified ``arithmetic'' of categorical idempotents.  To be more precise, we consider the linear algebra situation.  Let $R$ be a ring.  Some constructions which we would like to understand on the categorical level are:
\begin{enumerate}
\item If $e,f\in R$ are a pair of commuting idempotents, then $ef$ and $e+f-ef$ are idempotents.
\item if $e,f\in R$ are idempotents such that $fe =ef=e$, then $d(f,e) = f-e$ is an idempotent.
\item an idempotent decomposition of $1\in R$ is a collection of elements $e_i\in R$, indexed by a finite set $I$, such that
\begin{enumerate}
\item $e_ie_j = 0$ if $i\neq j$.
\item $e_i^2=e_i$ for all $i\in I$.
\item $1=\sum_i e_i$.
\end{enumerate}
Of course, (b) follows from (a) and (c).  If $J\subset I$ is any subset, then $e_J:= \sum_{i\in J} e_j$ is an idempotent, with complementary idempotent $e_{I\setminus J}$.
\end{enumerate}
In \S \ref{subsec:orderProps} and \S \ref{subsec:subquotients} we give the categorical analogues of (1) and (2) above, respectively.  Statement (3) also has a nice categorical analogue, except for one very important difference: any pair of complementary (categorical) idempotents comes equipped with a preferred order, since there are no nonzero maps from the counital idempotent to its unital complement.  Similarly, one expects the terms in a more general decomposition of $\one\in \AS$ into orthogonal idempotents to satisfy a similar semi-orthogonality condition with respect to morphisms.  The resulting partial order on idempotents is an intriguing phenomenon which is not visible on the decategorified level.

The categorical analogue of the completeness relation $1=\sum_i e_i$ is best understood in terms of the language of Postnikov systems, which we discuss in \S \ref{subsec:postnikov}.  We introduce the basics of categorical decompositions of identity in \S \ref{subsec:decompOfOne} and \S \ref{subsec:generalizedDecomp}.

\subsection{Suprema and infema}
\label{subsec:orderProps}
Let $(\unital_i,\eta_i)$ be unital idempotents in a triangulated monoidal category ($i=1,2$).  Recall Theorem \ref{thm:idemptOrder}, which states that $\unital_1\otimes \unital_2\cong \unital_1$ if and only if $\unital_2\otimes \unital_1\cong \unital_1$, which occurs if and only if there is a map $\phi:\unital_2\rightarrow \unital_1$ such that $\eta_1 = \phi\circ \eta_2$.  If either of these statements are true, we write $\unital_1\leq \unital_2$.  We call $\phi:\unital_2\rightarrow \unital_1$ is the canonical map.  Similar results hold for counital idempotents.  The purpose of this section is to prove:

\begin{proposition}\label{prop:suprema}
The partially ordered set of unital idempotents in $\AS$ has certain infema and suprema.  If $\{\unital_s\}_{s\in S}$ is a finite collection of commuting unital idempotents, then:
\begin{enumerate}
\item The tensor product $\bigotimes_s \unital_s$ has the structure of a unital idempotent; it is the infemum of $\{\unital_s\}_{s\in S}$.
\item For each $s\in S$, let $\unital_s^c$ be the complementary idempotent to $\unital_s$, and let $\unital$ be the complement of $\bigotimes_s \unital_s^c$.  Then $\unital=\sup\{\unital_s\}_{s\in S}$ is the supremum of $\{\unital_s\}_{s\in S}$.
\end{enumerate}
Similar statements hold for counital idempotents.
\end{proposition}
Now is a good time to recall Proposition \ref{prop:commuteWithidemp}, which implies that since the unital idempotents $\unital_s$ commute with one another, the counital idempotents $\unital_s^c$ commute with one another and with the $\unital_s$.

\begin{proof}
We first prove (2).  By an easy induction it suffices to assume that $S=\{1,2\}$.  That is, suppose we are given unital idempotents $\unital_1$ and $\unital_2$ which commute with each other.   Set $\counital:=\unital_1^c\otimes \unital_2^c$.  This object comes with a map $\e_1\otimes \e_2:\counital\rightarrow \one$, where $\e_i$ are the counits of $\unital_i^c$.  Set $\unital:=\Cone(\e_1\otimes \e_2)$.  Lemma \ref{lemma:coneOfTensor} below shows that $\unital$ fits into a distinguished triangle:
\[
\unital_1[-1]\rightarrow \unital_1^c\otimes \unital_2\rightarrow \unital \rightarrow \unital_1
\]
The first term above is annihilated on the left and right by $\unital_1^c$, and the second term is annihilated on the left and right by $\unital_2^c$.  Thus both terms are annihilated on the left and right by $\counital$.   It follows that $\unital$ is annihilated by $\counital$, which proves that $(\counital,\unital)$ are complementary idempotents.
 
We now show that $\unital_i\leq \unital$ for $i=1,2$, and $\unital$ is minimum among all idempotents with this property.

Since all the relevant complexes commute with one another, we have that $\unital^c=\unital_1^c \otimes \unital_2^c$ is annihilated by $\unital_1$ and $\unital_2$ on the left and right, hence tensoring with $\unital$ on the left or right fixes $\unital_1,\unital_2$.  That is to say, $\unital_1,\unital_2\leq \unital$.

Now, suppose $\idempotent$ is any unital idempotent such that $\unital_1,\unital_2\leq \idempotent$.  Let $\idempotent^c$ be the complementary idempotent to $\idempotent$.  By assumption, $\unital_i$ absorbs $\idempotent$ on the left and right for $i=1,2$.  Repeated applications of projector absorbing and its converse yield:
\begin{itemize}
\item $\idempotent^c\otimes \unital_i\cong 0$ ($i=1,2$).
\item $\idempotent^c\otimes \unital_i^c\cong \idemp^c$ ($i=1,2$), hence $\idemp^c\otimes \unital^c\cong \idemp^c$.
\item $\idempotent^c\otimes \unital\cong 0$.
\item $\idemp\otimes \unital \cong \idemp$ 
\end{itemize}
That is to say $\unital\leq \idempotent$, as desired.  A similar argument proves (1).
\end{proof}

 The following lemma was used in the proof of Proposition \ref{prop:suprema}:

\begin{lemma}\label{lemma:coneOfTensor}
Let $f:A\rightarrow B$ and $f':A'\rightarrow B'$ be morphisms in  $\AS$.  Then
\begin{eqnarray*}
\Cone(f\otimes g)
&\cong & \Cone\Big(B\otimes \Cone(f')[-1]\rightarrow  \Cone(f)\otimes A'\Big)\\
&\cong & \Cone\Big(\Cone(f)\otimes B'[-1]\rightarrow  A\otimes \Cone(f')\Big)
\end{eqnarray*}
\end{lemma}
\begin{proof}
Choose distinguished triangles
\[
\begin{diagram}[small]
A &\rTo^f & B & \rTo^g & C & \rTo^h & A[1]
\end{diagram}
\text{ \ \ \ \ \ \ and \ \ \ \ \ \ \ }
\begin{diagram}[small]
A' &\rTo^{f'} & B' & \rTo^{g'} & C' & \rTo^{h'} & A'[1]
\end{diagram}.
\]
Consider the diagram
\[
\begin{tikzpicture}
\tikzstyle{every node}=[font=\small]
\node (aa) at (0,0) {$A\otimes A'$};
\node (ba) at (3,0) {$B\otimes A'$};
\node (ca) at (6,0) {$C\otimes A'$};
\node (da) at (9,0) {$A\otimes A'[1]$};
\node (ab) at (0,-2.5) {$A\otimes A'$};
\node (bb) at (3,-2.5) {$B\otimes B'$};
\node (cb) at (6,-2.5) {$X$};
\node (db) at (9,-2.5) {$A\otimes A'[1]$};
\node (ac) at (0,-5) {$0$};
\node (bc) at (3,-5) {$B\otimes C'$};
\node (cc) at (6,-5) {$B\otimes C'$};
\node (dc) at (9,-5) {$0$};
\node (ad) at (0,-7.5) {$A\otimes A'[1]$};
\node (bd) at (3,-7.5) {$B\otimes A'[1]$};
\node (cd) at (6,-7.5) {$C\otimes A'[1]$};
\node (dd) at (9,-7.5) {$A\otimes A'[2]$}
;\\
\path[->,>=stealth',shorten >=1pt,auto,node distance=1.8cm, thick]
(aa) edge node {$f\otimes \Id$} (ba)
(ba) edge node {$g\otimes \Id$} (ca)
(ca) edge node {$h\otimes \Id$} (da)
(ab) edge node {$f\otimes g$} (bb)
(bb) edge node {} (cb)
(cb) edge node {} (db)
(ad) edge node {$f\otimes \Id[1]$} (bd)
(bd) edge node {$g\otimes \Id[1]$} (cd)
(cd) edge node {$h\otimes \Id[1]$} (dd)
(aa) edge node {$\Id$} (ab)
(ab) edge node {} (ac)
(ac) edge node {} (ad)
(ba) edge node {$\Id\otimes f'$} (bb)
(bb) edge node {$\Id\otimes g'$} (bc)
(bc) edge node {$\Id\otimes h'$} (bd)
(da) edge node {$\Id$} (db)
(db) edge node {} (dc)
(dc) edge node {} (dd)
;
\path[dashed,->,>=stealth',shorten >=1pt,auto,node distance=1.8cm, thick]
(ac) edge node {} (bc)
(bc) edge node {} (cc)
(cc) edge node {} (dc)
(ca) edge node {} (cb)
(cb) edge node {} (cc)
(cc) edge node {} (cd)
;
\end{tikzpicture}
\]
The top two rows a distinguished, as are the left two columns.  The $3\times 3$ lemma (Lemma 2.6 in \cite{MayTraces}) ensures the existence of the dashed arrows which make the third row and column distinguished, and such that the resulting diagram commutes (except for the bottom right square, which commutes up to sign).  

Distinguishedness of the second row means $X\cong \Cone(f\otimes g')$.  The third column then provides the desired distinguished triangle relating $\Cone(f\otimes f')$, $\Cone(f)\otimes A'$, and $B\otimes \Cone(f')$.  This proves the first statement.  The second is similar.
\end{proof}

\begin{definition}\label{def:supAndInf}
If $\unital_1$ and $\unital_2$ are commuting unital idempotents, then we write $\unital_1\wedge \unital_2 = \unital_1\otimes \unital_2 = \inf(\unital_1,\unital_2)$.  Let $\unital_1\vee\unital_2=\sup(\unital_1,\unital_2)$. 
\end{definition}
 Similar definitions apply to counital idempotents.  We have $(\unital_1\wedge \unital_2)^c \cong \unital_1^c\vee \unital_2^c$ and $(\unital_1\vee \unital_2)^c \cong \unital_1^c\wedge \unital_2^c$. 
 
We conclude this section with the following Mayer-Vietoris property:

\begin{proposition}\label{prop:MayerVietoris}
If $\UB$ and $\VB$ are commuting unital idempotents, then there is a distinguished triangle
\[
\begin{diagram}
\UB\vee \VB &\rTo^{\smMatrix{f\\g}} & \UB\oplus \VB &\rTo^{\smMatrix{h & -i}} & \UB\wedge \VB & \rTo &(\UB\vee\VB)[1],
\end{diagram}
\]
where the maps $f,g,h,i$ are the canonical morphisms of comparable unital idempotents (implied by Theorem \ref{thm:idemptOrder}).  Similarly, if $\CB$ and $\DB$ are counital idempotents, then there is a distinguished triangle
\[
\begin{diagram}
\CB \wedge \DB &\rTo^{\smMatrix{j & k}} & \CB\oplus \DB & \rTo^{\smMatrix{\ell \\ -m}} & \CB\wedge \DB & \rTo &(\CB\vee\DB)[1],
\end{diagram}
\]
where, again, the maps $j,k,\ell,m$ are the canonical morphisms of comparable counital idempotents.
\end{proposition}
\begin{proof}
This is Theorem 3.13 in \cite{BalmerFavi11}.  We warn that \emph{right idempotent} in \emph{loc.~cit.~} is what we call a unital idempotent; the partial order on right idempotents, is \emph{opposite} the partial order considered here.  Thus, what we call $\unital\vee \VB$ is called $\unital\wedge \VB$ in \emph{loc.~cit.~}.
\end{proof}

\subsection{Locally unital idempotents, excision}
\label{subsec:subquotients}
Starting in this section, we will focus our attention on unital idempotents in order to avoid stating every proposition twice.  The obvious modifications will produce the corresponding dual discussion concerning counital idempotents.

\begin{definition}\label{def:locallyUnital}
Let $\unital\leq \VB$ be unital idempotents in $\AS$.  Let $\nu:\VB\rightarrow \UB$ be the canonical map, and define an object $\DB(\VB,\UB)$ by choosing a distinguished triangle
\begin{equation}\label{eq:relativeTriangle}
\DB(\VB,\UB) \rightarrow \VB \buildrel \nu\over \rightarrow \UB \rightarrow \DB(\VB,\UB)[1]
\end{equation}
Then $\DB(\VB,\UB)$ is a counital idempotent relative to $\VB$.  We will refer to an idempotent obtained in this way as a \emph{locally unital idempotent} (called a locally closed idempotent in \cite{BoyDrin-idemp}).
\end{definition}

One might also call $\DB(\VB,\UB)$ a relatively counital idempotent in order to emphasize that it is a counital idempotent relative to a unital idempotent (in this case $\VB$).  In this section we prove a number of properties of such idempotents.  First we have a relative version of the fundamental theorem.

\begin{proposition}\label{prop:tripleOfRelIdempts}
If $\unital_1\leq \unital_2\leq \unital_3$, then there is a distinguished triangle
\begin{equation}\label{eq:Dtriang}
\DB(\unital_3,\unital_2)\rightarrow \DB(\unital_3,\unital_1)\rightarrow \DB(\unital_2,\unital_1)\rightarrow \DB(\unital_3,\unital_2)[1],
\end{equation}
and any two such triangles are canonically isomorphic in the following sense.  Suppose that there exists a distinguished triangle
\begin{equation}\label{eq:Etriang}
\EB_{32}\rightarrow \EB_{31} \rightarrow \EB_{21}\rightarrow \EB_{32}[1],
\end{equation}
such that $\DB(\unital_3,\unital_1)\cong \EB_{31}$, and choose an isomorphism $\psi:\DB(\unital_3,\unital_1)\rightarrow \EB_{31}$.  Then the following are equivalent:
\begin{enumerate}
\item $\DB(\unital_2,\unital_1)\otimes \EB_{21}\cong \DB(\unital_2,\unital_1)$.
\item $\EB_{21}\otimes \DB(\unital_2,\unital_1) \cong \DB(\unital_2,\unital_1)$.
\item $\DB(\unital_3,\unital_2)\otimes \EB_{32}\cong \EB_{32}$.
\item $\EB_{32}\otimes \DB(\unital_3,\unital_2)\cong \EB_{32}$.
\end{enumerate}
If either of these is satisfied then there is a unique morphism of triangles from (\ref{eq:Dtriang}) to (\ref{eq:Etriang}) which extends $\psi$.
\end{proposition}
\begin{proof}
Let $\nu_{ij}:\unital_j\rightarrow \unital_i$ denote the canonical maps ($i\leq j$).  Then $\nu_{13}=\nu_{12}\circ \nu_{23}$.  The octahedral axiom implies that $\DB(\unital_3,\unital_1)\cong \Cone(\nu_{13})$ fits into a distinguished triangle involving $\DB(\unital_2,\unital_1)\cong \Cone(\nu_{12})$ and $\DB(\unital_3,\unital_2)\cong \Cone(\nu_{23})$, as in the statement. 

The remaining statements follow by obvious modifications of the proof of Theorem \ref{thm:idemptOrder}.
\end{proof}

\begin{lemma}\label{lemma:subquotient}
$\DB(\VB,\UB)\cong \VB\otimes \UB^c$ for all unital idempotents $\UB\leq \VB$.
\end{lemma}
\begin{proof}
By the last statement of Theorem \ref{thm:idemptOrder}, $\nu$ realizes $\unital$ as a unital idempotent relative to $\VB$, hence $\DB(\VB,\UB)$ is the complementary idempotent.  In particular, $\DB(\VB,\UB)$ and $\UB$ are orthogonal with respect to $\otimes$.  Thus $\DB(\VB,\UB)\otimes \UB^c\cong \DB(\VB,\UB)$.  Now, tensor the defining triangle (\ref{eq:relativeTriangle}) on the right with $\UB^c$, obtaining a distinguished triangle $\DB(\VB,\UB)\rightarrow \VB\otimes \UB^c\rightarrow 0$.  The Lemma follows.
\end{proof}

\begin{proposition}[Generalized excision]\label{prop:excision}
If $\AB$ and $\BB$ are commuting unital idempotents, then $\DB(\AB\vee \BB,\AB)\cong \DB(\BB,\AB\wedge \BB)$.
\end{proposition}
\begin{proof}
By Proposition \ref{prop:MayerVietoris}, we have a distinguished triangle
\[
\AB\vee \BB\rightarrow \AB\oplus \BB\rightarrow \AB\wedge \BB\rightarrow \AB\vee \BB[1].
\]
Tensoring with $\AB^c$ and noting that $\AB$ and $\AB\wedge \BB = \AB\otimes \BB$ are annihilated by $\AB^c$, we obtain a distinguished triangle
\[
(\AB\vee \BB)\otimes \AB^c \rightarrow \BB\otimes \AB^c \rightarrow 0 \rightarrow (\AB\vee \BB)\otimes \AB^c[1].
\]
This implies $(\AB\vee \BB)\otimes \AB^c \cong \BB\otimes \AB^c$.  A similar argument shows $\BB\otimes(\BB^c\vee \AB^c)\cong \BB\otimes \AB^c$.  Thus, Lemma \ref{lemma:subquotient} gives us
\[
\DB(\AB\vee \BB, \AB) \cong (\AB\vee \BB)\otimes \AB^c \cong \BB\otimes \AB^c
\]
and
\[
\DB(\BB,\BB\wedge \AB) \cong \BB \otimes (\BB \wedge \AB)^c \cong \BB\otimes(\BB^c\vee \AB^c) \cong \BB\otimes \AB^c.
\]
This completes the proof.
\end{proof}

Now we show that every categorical idempotent in $\AS$ is isomorphic to $\DB(\VB,\UB)$ for some unital idempotents $\UB\leq \VB$.  First, a  lemma:

\begin{lemma}\label{lemma:relRelIdempts}
We have
\begin{enumerate}
\item Let $\eta:\one\rightarrow \unital$ be a unital idempotent and $\nu:\unital\rightarrow \idempotent$ be a unital idempotent relative to $\unital$.  Then $\nu\circ \eta:\one\rightarrow \idempotent$ is a unital idempotent.
\item   Let $\e:\counital\rightarrow \one$ be a counital idempotent and $\theta:\idempotent \rightarrow \counital$ be a counital idempotent relative to $\counital$.  Then $\e\circ \theta:\idempotent\rightarrow \one$ is a counital idempotent.
\end{enumerate}
\end{lemma}
\begin{proof}
We first prove (1).  Assume that $\idempotent\in \unital\AS\unital$, and $\nu:\unital\rightarrow \idempotent$ makes $\idempotent$ into a unital idempotent relative to $\unital$. The octahedral axiom implies that $\Cone(\nu\circ \eta)$ is isomorphic to the mapping cone on a map
\[
\Cone(\nu)[-1]\rightarrow \Cone(\eta).
\]
The second term is annihilated on the left and right by $\unital$ since $\eta:\one\rightarrow \unital$ is a unital idempotent, hence is also annihilated by $\idempotent$.  The first term is annihilated by $\idempotent$ since we are assuming that $\nu:\unital\rightarrow \idempotent$ is a unital idempotent rel $\unital$.  It follows that $\Cone(\nu\circ \eta)$ is annihilated on the left and right by $\idempotent$, hence $(\idempotent, \nu\circ \eta)$ is a unital idempotent.  This proves (1).  The proof of (2) is similar.
\end{proof}

\begin{proposition}\label{prop:classificationOfCatIdempts}.
Every categorical idempotent in $\AS$ isomorphic to $\DB(\VB,\UB)$ for some unital idempotents $\UB\leq \VB$. 
\end{proposition}
\begin{proof}
This is proven by an easy induction on depth.  The inductive step is provided by the following argument.  In the base case, note that every unital idempotent $\unital$ idempotent is isomorphic to $\DB(\unital,0)$, and the complementary counital idempotent is isomorphic $\DB(\one,\unital)$. 

The inductive step is taken care of by the following argument.  Let $\unital_1\leq \unital_3$ be unital idempotents, and suppose we are given some distinguished triangle
\[
\EB \rightarrow \DB(\unital_3,\unital_1)\rightarrow \FB\rightarrow \EB[1]
\]
with $\EB\otimes \FB\cong 0 \cong \FB \otimes \EB$.  That is, $\EB$ and $\FB$ are a unital and a counital idempotent with respect to $\DB(\unital_3,\unital_1)$.   We will construct a unital idempotent $\unital_2$ such that $\EB\cong \DB(\unital_3,\unital_2)$ and $\FB\cong \DB(\unital_2,\unital_1)$.

Since $\EB$ is a counital idempotent relative to $\DB(\unital_3,\unital_1)$, and $\DB(\unital_3,\unital_1)$ is a counital idempotent relative to $\unital_3$, Lemma \ref{lemma:relRelIdempts} says that $\EB$ has the structure of a counital idempotent relative to $\unital_3$.  Thus, we have the complementary unital idempotent $\unital_2$ relative to $\unital_3$ defined by choosing a distinguished triangle
\[
\EB\rightarrow \unital_3\rightarrow \unital_2\rightarrow \EB[1].
\]
Lemma \ref{lemma:relRelIdempts} says that $\unital_2$ is a unital idempotent since $\unital_3$ is.  Thus, $\EB\cong \DB(\unital_3,\unital_2)$.  Proposition \ref{prop:tripleOfRelIdempts} implies that the complementary idempotent (relative to $\DB(\unital_3,\unital_1)$ is $\FB\cong \DB(\unital_2,\unital_1)$, as claimed.
\end{proof}

Let's now take a moment to carefully describe the sense in which $\DB(\UB,\VB)$ is unique.  
\begin{proposition}\label{prop:relativeUniqueness}
Let $(\unital_i,\eta_i)$ and $(\VB_i,\nu_i)$ be unital idempotents with $\unital_i\leq \VB_i$ ($i=1,2$).  If, in addition, $\unital_1\leq \unital_2$ and $\VB_1\leq \VB_2$, then there is a unique map of distinguished triangles
\[
\begin{diagram}
\DB(\VB_2,\unital_2) & \rTo^{\e_2} & \VB_2 &\rTo^{\phi} & \unital_2 &\rTo & \DB(\VB_2,\unital_2)[1]\\
\dTo^{\psi} && \dTo^{\phi} && \dTo^{\phi} && \dTo^{\psi[1]}\\
\DB(\VB_1,\unital_1) & \rTo^{\e_1} & \VB_1 &\rTo^{\phi} & \unital_1 &\rTo & \DB(\VB_1,\unital_1)[1]
\end{diagram}
\]
in which the abused symbol $\phi$ denotes the canonical map associated to comparable unital idempotents.   If $\DB(\VB_1,\unital_1)\cong \DB(\VB_2,\UB_2)$, then $\psi$ is an isomorphism.
\end{proposition}
\begin{proof}
Let us abbreviate $\DB_i:=\DB(\VB_i,\UB_i)$ throughout the proof.  When we want to be explicit, we will let $\phi_{\UB,\VB}:\VB\rightarrow \UB$ denote the canonical map associated to idempotents $\UB\leq \VB$.  The existence of a map of triangles follows from the axioms of triangulated categories (the square involving the canonical maps $\phi$ commutes).  The canonical maps $\phi$ are unique.  Now we argue that $\psi$ is unique.  We first claim that, under the hypotheses, $\DB_i, \VB_i, $ and $\UB_i$ are all fixed by tensoring with $\VB_2$ on the left or right.  Indeed, $\VB_2$ is the maximum among $\{\unital_i,\VB_i\}$, hence $\UB_i,\VB_i\in \VB_2\AC\VB_2$.  This category is triangulated, hence $\DB_i\in \VB_2\AC\VB_2$ as well.  Thus, without loss of generality we may assume that $\VB_2=\one$.

Now $\DB_2$ is a counital idempotent with complement $\unital_2$.  Since $\unital_2\geq \unital_1$, $\DB_2\otimes \unital_2\cong 0$ implies $\DB_2\otimes \unital_1\cong 0$ as well.  Thus, Corollary \ref{cor:unitalCounitalOrtho} implies that $\Homg(\DB_2,\unital_1)\cong 0$.  Applying $\Homg(\DB_2,-)$ to the distinguished triangle $\DB_1\rightarrow \VB_1\rightarrow \UB_1$ gives a long exact sequence in which every third term is zero.  Thus $\e_1$ induces an isomorphism $\Homg(\DB_2,\DB_1)\rightarrow \Homg(\DB_2,\VB_1)$.  In particular, there is a unique $\psi:\DB_2\rightarrow \DB_1$ such that $\e_1\circ \psi = \phi_{\VB_1,\VB_2}\circ \e_2$.  This proves uniqueness of $\psi$.

Now, suppose that $\DB_1 \cong \DB_2$.  We must show that $\Cone(\psi)\cong 0$.   First, consider $\Cone(\e_1\circ \psi)$.  By the octahedral axiom this fits into a distinguished triangle of the form
\[
\Cone(\psi)\rightarrow \Cone(\e\circ \psi)\rightarrow \Cone(\e)\rightarrow \Cone(\psi)[1].
\]
Note that $\DB_2\otimes \Cone(\e_1)\cong \DB_2\otimes \unital_1\cong 0$.  It follows that
\[
\DB_2\otimes \Cone(\psi) \cong \DB_2\otimes \Cone(\e\circ \psi) = \DB_2\otimes  \Cone(\phi_{\VB_1,\VB_2}\circ \e_2).
\]
On the other hand note that $\psi$ is a map from $\DB_2$ to $\DB_1\cong \DB_2$, hence $\Cone(\psi)\in \DB_2\AS$.  Thus, the left most object above is isomorphic to $\Cone(\psi)$.  To show that $\psi$ is an isomorphism, we must show that the right-most term above is zero. 

By the octahedral axiom, $\Cone(\phi_{\VB_1,\VB_2}\circ \e_2)$ fits into a distinguished triangle
\[
\Cone(\e_2)\rightarrow \Cone(\phi_{\VB_1,\VB_2}\circ \e_2) \rightarrow \Cone(\phi_{\VB_1,\VB_2})\rightarrow \Cone(\e_2)[1].
\]
Once again, $\DB_2\otimes \Cone(\e_2)\cong \DB_2\otimes \UB_2\cong 0$, hence tensoring this distinguished triangle with $\DB_2$ gives
\[
\DB_2\otimes \Cone(\phi_{\VB_1,\VB_2}\circ \e_2) \cong \DB_2\otimes \Cone(\phi_{\VB_1,\VB_2}) .
\]
Now, $\DB_1\cong \DB_2$ satisfies $\DB_2\otimes \VB_1 \cong \DB_2$.  But $\VB_1$ is orthogonal to $\Cone(\phi_{\VB_1,\VB_2})$, since the two are orthogonal idempotents relative to $\VB_2$.  It follows that $\DB_2$ is orthogonal to $\Cone(\phi_{\VB_1,\VB_2})$, which proves that $\Cone(\psi)\cong 0$.  This completes the proof.
\end{proof}

\subsection{Postnikov systems}
\label{subsec:postnikov}
Let $0=\unital_0\leq \unital_1\leq \cdots \leq \unital_n$ be a sequence of unital idempotents in $\AS$, and let $\idemp_i:=\DB(\unital_i,\unital_{i-1})$ for all $1\leq i\leq n$.  Then repeated application of Proposition \ref{prop:tripleOfRelIdempts} says that $\unital_n$ can be written as a cone involving $\idemp_1$ and $\DB(\unital_n,\unital_{1})$, and that this latter term can be written as a cone involving $\idemp_2$ and $\DB(\unital_n,\unital_2)$, and so on.

Such sequences of iterated mapping cones are called Postnikov systems, which we now discuss.  Consider the following situation: let $\TS$ be a triangulated category, and let $F_n\rightarrow F_{n-1}\rightarrow \cdots \rightarrow F_0$ be a directed system of objects in $\TS$; let us denote the maps by $\sigma_i:F_i\rightarrow F_{i-1}$.  Consider the following construction: for each $0\leq i\leq n$, define an object $R_i$ by choosing a distinguished triangle
\begin{equation}\label{eq:convTriang}
F_{i+1}\buildrel \sigma_{i+1}\over \longrightarrow F_i \buildrel \pi_i\over \rightarrow R_i\buildrel \Delta_i\over \longrightarrow F_{i+1}[1].
\end{equation}
By convention, $F_{n+1}=0$, hence we may choose $R_n=F_n$ and $\pi_n=\Id$.  Define $d_i=\pi_{i+1}[1]\circ \Delta_i : R_i\rightarrow R_{i+1}[1]$.  The composition of successive $d$'s is zero since
\[
d_{i+1}[1]\circ d_i = \pi_{i+2}[2]\circ \Delta_{i+1}[1]\circ \pi_{i+1}[1]\circ \Delta_i
\]
The middle two terms contribute $(\Delta_{i+1}\circ \pi_{i+1})[1]$ to the right-hand-side; this map is zero since the composition of successive maps in a distinguished triangle is zero.  Thus, we have a chain complex
\begin{equation}\label{eq:Rcomplex}
R_0\rightarrow R_1[1]\rightarrow\cdots \rightarrow  \cdots R_n[n].
\end{equation}
This is to be thought of as an object of $\KC(\TS)$. 

There is a related construction, which goes as follows.  Retain notation as above.  For each $0\leq i\leq n$, define an object $L_i$ by choosing a distinguished triangle
\begin{equation}\label{eq:Lcomplex}
L_i \longrightarrow F_i \buildrel \sigma_i \over \longrightarrow F_{i-1} \longrightarrow L_i[1]
\end{equation}
By convention, $F_{-1}=0$ and $L_0=F_0$.  A similar argument to the above produces a chain complex $L_0\rightarrow L_1\rightarrow \cdots L_n[n]$.  By uniqueness of mapping cones, and rotation of triangles, the objects $L_i$ are related to the $R_i$ by $L_i=R_{i-1}[1]$ for $1\leq i\leq n-1$, in addition to $L_0=F_0$ and $R_n=F_n$.  We will say that $F_0$ is a \emph{right convolution} of the complex (\ref{eq:Rcomplex}), and $F_0$ is a \emph{left convolution} of the complex (\ref{eq:Lcomplex}).  

In this way, directed systems in $\TS$ give rise to complexes over $\TS$, in two ways.  The notion of Postnikov systems and convolutions captures the reverse processes.  Below we state a more precise definition:

\begin{definition}\label{def:shiftFriendly}
A \emph{shift-friendly complex} over $\TS$ is a collection of objects $M_0,M_1,\ldots,M_n\in \TS$ and maps $d_i:M_i\rightarrow M_{i+1}[1]$ such that $d_{i+1}[1]\circ d_i=0$.  We will indicate a shift-friendly complex by the notation
\begin{equation}\label{eq:shiftFriendlyM}
M_0\buildrel[1]\over \rightarrow M_1 \buildrel[1]\over \rightarrow \cdots \buildrel[1]\over \rightarrow M_n
\end{equation}
\end{definition}
Note that shift-friendly complex $M_\bullet$ gives rise to an honest complex $M_\bullet'\in \TS$ by setting $M_i':=M_i[i]$.

\begin{definition}\label{def:PostnikovSystem}
Let $M_\bullet$ be a shift-friendly complex over $\TS$.  A \emph{right Postnikov system attached to $M_\bullet$} is a collection of objects $F_i$ ($0\leq i\leq n$) with $F_n=M_n$, and a collection of distinguished triangles
\[
F_{i+1}\buildrel \sigma_{i+1}\over \longrightarrow F_i \buildrel \pi_i\over \rightarrow M_i\buildrel \Delta_i\over \longrightarrow F_{i+1}[1].
\]
such that the differential of $M_\bullet$ satisfies $d_i = \pi_{i+1}[1]\circ \Delta_i: M_i\rightarrow M_{i+1}$ for all $i$.   By convention, $F_{n+1}=0$, and we require that $\pi_n$ be the identity of $F_n=M_n$.   In this case we call $F_0$ a \emph{right convolution} of $M_\bullet$. 

A \emph{left Postnikov system is of $M_\bullet$} is defined similarly, except that we require $F_0=M_0$, and the distinguished triangles relating the $F_i$ and $M_i$ look like
\[
M_i\rightarrow F_i\rightarrow F_{i-1}\rightarrow M_i.
\]
In this case, $F_n$ is called a \emph{left convolution of $M_\bullet$}.  The collection of left convolutions and right convolutions of a given $M_\bullet$ are equal; we denote this set by $\Tot(M_\bullet)$.  This set can be empty.  Even though convolutions are not unique, we will sometimes write $F=\Tot(M_\bullet)$ to indicate that $F$ is a convolution of $M_\bullet$.
\end{definition}
\begin{remark}\label{rmk:shiftFriendly}
In most references, one finds non shift-friendly versions of the above definition.  If $M_\bullet'$ is an honest complex (not shift-friendly) and $M_i=M_i[-i]'$, then any convolution $M_\bullet'$ in the usual sense will also be a convolution of $M_\bullet$ in our, shift-friendly sense.  We will consider shift-friendly convolutions almost exclusively.  When we wish to emphasize that $F$ is a shift-friendly convolution, we will write
\[
F = \Tot(M_0\buildrel[1]\over \rightarrow M_1 \buildrel[1]\over \rightarrow \cdots \buildrel[1]\over \rightarrow M_n)
\]
Note that on the level of Grothendieck groups, we have $[F]=[M_0]+[M_1]+\cdots +[M_n]$.  In the non-shift-friendly version, we would have an alternating sum: $[\Tot(N_\bullet)]=\sum_{i=0}^n(-1)[N_i]$.
\end{remark}

\begin{example}\label{ex:convCx}
If $\TS\subset \KC(\CC)$ for some additive category $\CC$, we can give a more explicit description of convolutions.  Now, let $M_\bullet' =M_{0}'\rightarrow M_{1}'\rightarrow \cdots \rightarrow  M_n'$  be an object of $\KC(\KC(\CC))$ (non-shift friendly).  That is, $(M_\bullet,d)$ is a complex of complexes.  We require that the composition of successive differentials in $M_\bullet$ be homotopic to zero, rather than zero on the nose.  If $d_{i+1}\circ d_i$ were zero on the nose, then $M_\bullet$ would be a bicomplex, and we could form the total complex $\Tot(M_\bullet)$, which would be a convolution of $M_\bullet$ (hence the notation).

Let us shift now, to obtain the shift-friendly version of $M_\bullet'$.  That is, set $M_i=M_i'[-i]$.  A convolution of $M_\bullet$ is a choice of linear maps  $d_{i,j}\in \Hom^{1}_{\CC}(M_j,M_i)$ which are homogeneous of degree 1, ($i\geq j$) such that
\begin{itemize}
\item $d_{ii}=d_{M_i}$, the differential internal to $M_i$.
\item $\sum_{i\geq j} d_{i,j}$ is a differential on $\bigoplus_i M_i$.
\end{itemize}
To compare with what was said for general triangulated categories, let $F_0$ denote $\bigoplus_{i=M}^N M_i$ with its differential $\sum_{i\geq j} d_{i,j}$.  In general, let $F_k$ denote the subcomplex
\[
F_k=\bigoplus_{i=k}^n M_i\subset F_0
\]
Note that the differential of $F_k$ is not the boring differential, but rather it is $\sum_{i\geq j\geq k} d_{i,j}$.  Then $F_{i+1}$ is a subcomplex of $F_i$ with quotient $F_i/F_{i+1}\cong M_{i+1}[$.  We leave it to the reader to verify that this gives rise to a right Postnikov system associated to $M_\bullet$.  Similarly, the complexes $E_k:=\bigoplus_{i\leq k}M_i$ with differential $\sum_{k\geq i\geq j} d_{ij}$ form a left Postnikov system associated to $M_\bullet$.  
\end{example}

\begin{example}\label{ex:filtrationConv}
Let $A\in \CC$ be an object of an abelian category, and let $F_n\subset F_{n-1}\subset \cdots \subset F_0 =A$ be a descending filtration on $A$.  For each $i$, we have the subquotient $R_i = F_i/F_{i+1}$.  By convention, $F_{n+1}=0$ and $R_n=F_n$.  Then in the derived category $D(\CC)$ we have distinguished triangles as in (\ref{eq:convTriang}), hence $A$ is isomorphic in $D(\CC)$ to a convolution:
\[
A \cong \Tot(R_0 \buildrel[1]\over \rightarrow R_1 \buildrel[1]\over \rightarrow \cdots \buildrel[1]\over \rightarrow R_n)
\]
\end{example}

We have two useful constructions involving Postnikov systems. Below, we will prefer left Postnikov systems; similar considerations also apply to right Postnikov systems.

\begin{proposition}[Truncation]\label{prop:truncation}
Let $M_\bullet$ be a shift-friendly complex of the form (\ref{eq:shiftFriendlyM}) and $F_n\rightarrow F_{n-1}\rightarrow \cdots \rightarrow F_0$ a left Postnikov system attached to $M_\bullet$.  Denote the maps in this directed system by $f_{ij}:F_j\rightarrow F_i$, for $i\leq j$.  Then for each pair of integers $0\leq k\leq \ell\leq n$, we have a convolution of the truncated complex
\begin{equation}\label{eq:truncation}
\Cone(f_{k-1,\ell}) \cong \Tot\Big(M_k\buildrel[1]\over \rightarrow M_{k+1} \buildrel[1]\over \rightarrow \cdots \buildrel[1]\over \rightarrow M_\ell\Big).
\end{equation}
In fact, there is a left Postnikov system attached to the truncated shift-friendly complex with terms given by $G_i := \Cone(f_{k-1,i})$.
\end{proposition}
\begin{proof}
Let us define the objects $G_i$ ($k\leq i\leq \ell$) by choosing distinguished triangles
\[
G_i \buildrel g_i\over  \longrightarrow F_i\buildrel f_{k-1,i}\over  \longrightarrow F_{k-1} \longrightarrow G_i[1].
\]
We also set $G_{k-1}:=0$.  The axiom (T3) for triangulated categories, applied to the composition $f_{k-1,i}[-1]=f_{k-1,i-1}[-1]\circ f_{i-1,i}[-1]$ yields a distinguished triangle
\[
M_i \longrightarrow G_i \buildrel h_i\over \longrightarrow G_{i-1} \buildrel \d_i\over \longrightarrow  M_i[1].
\]
Examination of the commutative diagram associated to (T3) shows that the composition $M_i\rightarrow G_i\rightarrow M_{i+1}[1]$ agrees with the composition $M_i\rightarrow F_i\rightarrow M_{i+1}[1]$.  Thus, $G_\ell \rightarrow G_{\ell-1}'\rightarrow \cdots \rightarrow G_k$ is a left Postnikov system for the truncation $M_k\rightarrow \cdots \rightarrow M_\ell$, as claimed.
\end{proof}
The construction in the above proof  implies the existence of a distinguished triangle:
\begin{equation}\label{eq:truncTriang}
\Tot\Big(M_k\buildrel[1]\over \rightarrow \cdots \buildrel[1]\over \rightarrow M_\ell\Big) \rightarrow \Tot\Big(M_1\buildrel[1]\over \rightarrow \cdots \buildrel[1]\over \rightarrow M_\ell\Big)\rightarrow \Tot\Big(M_1\buildrel[1]\over \rightarrow \cdots \buildrel[1]\over \rightarrow M_{k-1}\Big) \buildrel[1]\over \rightarrow
\end{equation}

\begin{proposition}[Reassociation]\label{prop:reassociation}
Let $M_\bullet$ be a shift-friendly complex as in (\ref{eq:shiftFriendlyM}), and let $F_n\rightarrow F_{n-1}\rightarrow \cdots \rightarrow F_0$ be a left Postnikov system associated to $M_\bullet$.  Let $0=i_0< i_1\cdots < i_r=n$ be integers, and let
\[
N_k = \Tot\Big(M_{i_{k-1}+1} \buildrel [1]\over \rightarrow \cdots \buildrel [1]\over \rightarrow M_{i_k}\Big)
\]
be the truncated convolution constructed in Proposition \ref{prop:truncation} ($1\leq k\leq r$).  Then the $N_k$ form a shift-friendly complex, and
\[
F_0=\Tot(M_\bullet) \cong \Tot\Big(N_1  \buildrel [1] \over \rightarrow N_2 \buildrel [1] \over \rightarrow \cdots \buildrel [1] \over \rightarrow N_r\Big).
\]
\end{proposition}
\begin{proof}
For brevity, if $J\subset \{0,1,\ldots,n\}$ is an interval, then let $M_J$ denote the corresponding truncation of $M_\bullet$, and let $G_J=\Tot(M_J)$ denote the truncation constructed in Proposition \ref{prop:truncation}.  For each $1\leq k\leq r$, let $J_k = \{i_{k-1}+1,\ldots,i_k\}$, so that $\{0,1,\ldots,n\}=J_1\sqcup \cdots \sqcup J_r$ and $N_k =\Tot(M_{J_k})) = G_{J_k}$.

Consider the directed system $G_{J_1\sqcup \cdots \sqcup J_r} \rightarrow   G_{J_1\sqcup \cdots \sqcup J_{r-1}} \rightarrow \cdots \rightarrow  G_{J_1}$.  From the construction in Proposition \ref{prop:truncation} (see (\ref{eq:truncTriang}), these fit into distinguished triangles
\[
G_{J_s} \rightarrow G_{J_1\sqcup \cdots \sqcup J_s} \rightarrow G_{J_1\sqcup \cdots \sqcup J_{s-1}} \rightarrow G_{J_s}[1],
\]
Since $G_{J_s}=N_s$, we have constructed the desired Postnikov system.
\end{proof}

\subsection{Linear decompositions of identity}
\label{subsec:decompOfOne}

\begin{definition}\label{def:linearDecomp}
A (linear) \emph{idempotent decomposition of identity} is a sequence of unital idempotents $\{\unital_i,\eta_i\}_{i=0}^n$ such that
\begin{enumerate}
\item $i\leq j$ implies $\unital_i\leq \unital_j$
\item $\unital_n=\one$.  
\end{enumerate}
The idempotents $\idemp_i:=\DB(\unital_i,\unital_{i-1})$ will be called the \emph{atomic subquotients} of $\{\unital_i\}_{i=0}^n$.  By convention, $\unital_{-1}=0$ and $\idemp_0=\unital_0$.
\end{definition}
The canonical maps give rise to a directed system $\one \cong \unital_n\rightarrow \unital_{n-1}\rightarrow \cdots \rightarrow \unital_0$. Given the distinguished triangles $\DB(\unital_i,\unital_{i-1})\rightarrow \unital_i\rightarrow \unital_{i-1}\rightarrow \DB(\unital_i,\unital_{i-1})[1]$, it follows that there is a left Postnikov system with convolution $\one$:
\begin{equation}\label{eq:resOfId}
\one \cong \Tot\Big(\idemp_0 \buildrel [1]\over \rightarrow \idemp_1\buildrel [1]\over \rightarrow \cdots \buildrel [1]\over \rightarrow \idemp_n\Big).
\end{equation}
Similarly, the complementary idempotents form a right Postnikov system $\unital_{n-1}^c \rightarrow \cdots \unital_0^c \rightarrow \unital_{-1}^c\cong \one $ with convolution $\one$.

Note that on the level of Grothendieck group we have $[\one]=\sum_{i=0}^n [\idempotent_i]$.
\begin{proposition}\label{prop:subquotientProps}
Let $\{\unital_i,\eta_i\}_{i=0}^n$ be an idempotent decomposition of identity with subquotients $\idemp_i$.  Then:
\begin{enumerate}
\item the $\idempotent_i$ are mutually orthogonal idempotents.
\item $\Homg(\idemp_i,\idemp_j)=0$ unless $i\leq j$.
\end{enumerate}
\end{proposition}
\begin{proof}
By construction, $\idemp_i$ absorbs $\unital_i$ and is killed by $\unital_{i-1}$ on the left and right.  Further, since $j\geq i$ implies $\unital_j\geq \unital_i$, it follows that $\idemp_i$ absorbs $\unital_j$ on the left and right for all $j\geq i$.   Thus, if $i<j$, then  the functors $\idemp_{j-1}\otimes(-)$ and $(-)\otimes \idemp_{j-1}$ fix $\idemp_i$ and annihilate $\idemp_j$.  Statement (1) follows.

For statement (2), let $i<j$. note that $\idemp_j$ is annihilated by $\unital_{j-1}$, hence is fixed by the complement $\unital_{j-1}^c$.  On the other hand $\idemp_i$ is fixed by $\unital_{j-1}$.  Thus (2) follows from the semi-orthogonality of idempotents (Theorem \ref{thm-semiOrtho}).
\end{proof}

Now, we want to show that the entire decomposition of identity is determined, up to canonical isomorphism, by the isomorphism classes of the $\idemp_i$.

\begin{theorem}\label{thm:linearDecomp}
Suppose $\AS$ is a triangulated monoidal category and $\idemp_0,\idemp_1,\ldots,\idemp_n\in \AS$ objects such that
\begin{enumerate}
\item $\idemp_i\otimes \idemp_j \cong 0$ for $i\neq j$, and
\item There exist maps $\idemp_i\rightarrow \idemp_{i+1}[1]$ which make $\idemp_\bullet$ into a shift-friendly complex (Definition \ref{def:shiftFriendly}) such that
\begin{equation}\label{eq:decompOfOne}
\one\cong \Tot(\idemp_0\buildrel[1]\over \rightarrow \cdots \buildrel[1]\over\rightarrow \idemp_n).
\end{equation}
\end{enumerate}
Then there is an idempotent decomposition of identity $\{\unital_i,\eta_i\}_{i=0}^n$ with $\DB(\unital_i,\unital_{i-1})\cong \idemp_i$.  Furthermore $\{\unital_i,\eta_i\}_{i=0}^n$ is determined by the isomorphism classes of the $\idemp_i$ up to canonical isomorphism.
\end{theorem}
\begin{proof}
Suppose there is some shift-friendly complex
\[
\idemp_\bullet = \idemp_0\buildrel[1]\over \rightarrow \cdots \buildrel[1]\over\rightarrow \idemp_n
\]
and a left Postnikov system $\unital_n\rightarrow \unital_{n-1}\rightarrow \cdots \rightarrow \unital_0$ attached to $\idemp_\bullet$ such that $\unital_n\cong \one$ as in the hypotheses.   Note that
\[
\unital_i = \Tot(\idemp_0\buildrel [1]\over \rightarrow \cdots \buildrel [1]\over \rightarrow \idemp_i).
\]
Let $\unital_i^c$ denote the truncation
\[
\unital_i^c = \Tot(\idemp_{i+1} \buildrel [1]\over \rightarrow \cdots \buildrel [1]\over \rightarrow \idemp_n).
\]
By construction, $\unital_i$ and $\unital_i^c$ are related by a distinguished triangle
\[
\unital_i^c\rightarrow \one \rightarrow \unital_i\rightarrow \unital_i^c[1].
\]
Now, since $\idemp_i\otimes \idemp_j\cong 0$ for $i\neq j$, it follows that $\unital_i$ and $\unital_i^c$ are orthogonal (see Lemma \ref{lemma:tensorConv}), hence they are a pair of complementary idempotents.  In fact $\unital_i$ is orthogonal to $\unital_j^c$ for every $j>i$.  This implies that $\unital_i\leq \unital_j$ for $i\leq j$.  Choosing an isomorphism $\one\cong \unital_n$ then gives a decomposition of identity $\{\unital_i,\eta_i\}$ as in the statement.

The question of uniqueness is addressed in Proposition \ref{prop:subquotientCharacterization} below.
\end{proof}

\begin{proposition}\label{prop:subquotientCharacterization}
Let $\unital_0\leq \unital_1\leq \cdots \leq \unital_n$ be unital idempotents and set $\idemp_i:=\DB(\unital_i,\unital_{i-1})$ for $1\leq i\leq n$, and $\idemp_0:=\unital_0$.  Then $\DB(\unital_j,\unital_i)$ is uniquely characterized up to canonical isomorphism by
\begin{enumerate}
\item $\DB(\unital_j,\unital_i)$ is in the triangulated hull $\langle\idemp_k \rangle_{i+1\leq k\leq j}$.
\item $\DB(\unital_j,\unital_i)\otimes(-)$ and $(-)\otimes \DB(\unital_j,\unital_i)$ fix up to isomorphism every object of $\langle\idemp_k \rangle_{i+1\leq k\leq j}$.
\end{enumerate}
\end{proposition}
\begin{proof}
We regard $\unital_0\leq \cdots \leq \unital_n$ as a left Postnikov system, where the maps $\unital_j\rightarrow \unital_i$ ($0\leq i\leq j\leq n$) are the canonical maps of comparable unital idempotents (Theorem \ref{thm:idemptOrder}).  Thus, we regard $\unital_n$ as a convolution:
\[
\unital_n \cong \Tot\Big(\DB(\unital_0,0)\buildrel[1]\over \longrightarrow \DB(\unital_1,\unital_0)\buildrel[1]\over \longrightarrow  \cdots \buildrel[1]\over \longrightarrow \DB(\unital_{n-1},\unital_n)\Big).
\]
That the terms are isomorphic to $\DB(\unital_i,\unital_{i-1})$ follows from the existence of distinguished triangles $\DB(\unital_i,\unital_{i-1})\rightarrow \unital_i\rightarrow \unital_{i-1}\rightarrow \DB(\unital_i,\unital_{i-1})[1]$.

Set $\idemp_i:=\DB(\unital_i,\unital_{i-1})$.  By Proposition \ref{prop:truncation}, the truncations of this Postnikov system are the cones of the given maps $\unital_j\rightarrow \unital_i$, with $0\leq i\leq j\leq n$.  But, again, the cone of the canonical map is isomorphic to $\DB(\unital_j,\unital_i)$ by definition.  We conclude that
\[
\DB(\unital_j,\unital_i)\cong \Tot\Big(\idemp_{i+1}\buildrel[1]\over \longrightarrow \cdots \buildrel[1]\over \longrightarrow \idemp_j \Big).
\]
Thus, $\DB(\unital_j,\unital_i)$ is in the triangulated hull $\langle\idemp_k \rangle_{i+1\leq k\leq j}$.  

By Proposition \ref{prop:reassociation} this can be reassociated into a convolution
\begin{equation}\label{eq:reassDecomp}
\one \cong \Tot\Big(\unital_i \buildrel[1]\over \longrightarrow \DB(\unital_j,\unital_i)\buildrel[1]\over \longrightarrow \DB(\unital_n,\unital_j) \Big),
\end{equation}
where we are using the fact that $\DB(\unital_i,0)\cong \unital_i$.  Note that if $i+1\leq k\leq j$, then $\idemp_k$ is annihilated by $\unital_i$ and $\unital_j^c$ by Lemma \ref{lemma:tensorConv}.  Thus, tensoring the convolution (\ref{eq:reassDecomp}) with $\idemp_k$ yields
\[
\DB(\unital_j,\unital_i)\otimes \idemp_k \cong \idemp_k\cong \idemp_k\otimes \DB(\unital_j,\unital_i).
\]
By Proposition \ref{prop-coneClosure}, the full subcategory consisting of objects fixed by $\DB(\unital_j,\unital_i)\otimes (-)$ or $(-)\otimes \DB(\unital_j,\unital_i)$ is triangulated.  Thus $\DB(\unital_j,\unital_i)$ satisfies properties (1) and (2) of the statement.

Now assume that $\FB$ satisfies properties (1) and (2).  Then $\FB\otimes \DB(\unital_j,\unital_i)$ is isomorphic to $\FB$ and $\DB(\unital_j,\unital_i)$, which proves uniqueness. 
\end{proof}

The following lemma was used repeatedly above; its proof is straightforward.
\begin{lemma}\label{lemma:tensorConv}
Let $\TS$ be a triangulated monoidal category, let $M_\bullet = M_0\rightarrow \cdots \rightarrow M_n$ be a chain complex over $\TS$, and let $N=\Tot(M_\bullet)$ be a convolution.  Then for any object $A\in \TS$, tensoring with $A$ gives a chain complex
\[
A\otimes M_\bullet = A\otimes M_0\buildrel \Id_A\otimes d\over\longrightarrow \cdots  \buildrel \Id_A\otimes d\over\longrightarrow  A\otimes M_n,
\]
with convolution $A\otimes N = \Tot(A\otimes M_\bullet)$.  In particular, if $A\otimes M_i\cong 0$ for all $i$, then $A\otimes N\cong 0$.  Similar statements hold with $A\otimes(-)$ replaced by $(-)\otimes A$.\qed
\end{lemma}

\subsection{Generalized decompositions of identity}
\label{subsec:generalizedDecomp}
In this section we consider idempotent decompositions of identity indexed by more interesting partially ordered sets.  But first we need some elementary combinatorial notions.  Let $(I,\leq)$ be a finite partially ordered set.  A subset $J\subset I$ is \emph{convex} if $j_1,j_2\in J$ and $j_1\leq i\leq j_2$ implies $i\in J$.  A subset $J\subset I$ is an \emph{ideal} (resp.~\emph{coideal}) if $j\in J$ and $i\leq j$ (resp.~$i\geq j$) implies $i\in J$.  Note that $J$ is a coideal iff $I\setminus J$ is an ideal, and any convex set can be written as the intersection of an ideal with a coideal.  The set of ideals (or coideals) $I$ is closed under unions and intersections, and forms a partially ordered set by inclusion. 

\begin{remark}
The set $I$ can be given a topology by declaring a subset $U\subset I$ to be open iff it is an ideal (resp.~a coideal).
\end{remark}

\begin{definition}\label{def:decompOfOne}
Let $(I,\leq)$ be a partially ordered set and $\AS$ a triangulated monoidal category.  An \emph{idempotent decomposition of identity} (indexed by $I$) is a collection of unital idempotents $\{\unital_J,\eta_J\}$ indexed by ideals $J\subset I$, such that
\begin{enumerate}
\item $\unital_I\cong \one$ and $\unital_\emptyset \cong 0$,
\item $\unital_{J\cap K}\cong \unital_J\otimes \unital_K$,
\item $\unital_{J\cup K} \cong \unital_J\vee \unital_K$ (Definition \ref{def:supAndInf}).
\end{enumerate}
Note that if $K\subset J\subset I$ are ideals, then $\unital_K\leq \unital_J$ by property (2).

Let $J\subset I$ be any convex set.  We can express $J$ as $J=K_2\cap K_1^c$ for some ideals $K_1\subset K_2$.  Define $\idemp_J=\DB(\unital_{K_2},\unital_{K_1})$.  This idempotent depends only on $J$ by excision (Proposition \ref{prop:excision}).  The idempotents $\idemp_J$ are called the subquotients of the decomposition of identity; if $J=\{i\}$ is a singleton, then $\idemp_i=\idemp_{\{i\}}$ will be called the \emph{atomic subquotients} of the decomposition of identity.
\end{definition}
The properties (1) and (2) imply that in order to specify an idempotent decomposition of identity, it is enough to define $\unital_J$ for \emph{basis} of $I$.  Here, a basis is a collection of ideals which generate the entire collection under the operations of $\cup$ and $\cap$.

\begin{remark}\label{rmk:sheavesOfIdempts}
If $X$ is a topological space, and $\TS$ is a triangulated category, then one has the notion of a \emph{sheaf of objects of $\TS$}, which is an object $\FS_U\in \TS$ for each open set $U\subset X$, together with a coherent family of restriction maps $\rho_{V,U}:\FS_U\rightarrow \FS_V$ whenever $V\subset U$.  The usual sheaf condition gets replaced by the requirement that $\FS_{U\cup V}$ be isomorphic to a (homotopy) pullback associated to the restriction maps $\FS_U,\FS_V\rightarrow \FS_{U\cap V}$.  

Thus, a generalized idempotent decomposition of identity in a triangulated monoidal category $\AS$ may be thought of as a \emph{sheaf of unital idempotents} over $I$, with its topology given by ideals.  The ``sheaf condition'' is guaranteed by property (3) of Definition \ref{def:decompOfOne} together with the Mayer-Vietoris property (Proposition \ref{prop:MayerVietoris}).
\end{remark}

\begin{example}\label{ex:generalizedConvolutionCx}
Suppose $\AS\subset \KC(\CC)$ for some additive monoidal category $\CC$.  Let $\idemp_i$ be chain complexes and suppose we are given degree 1 linear maps $d_{ij}\in \Hom^1_\CC(\idemp_j,\idemp_i)$ ($i\geq j$) such that
\begin{enumerate}
\item $d_{ii}=d_{\idempotent_i}$, the differential on $\idempotent_i$.
\item $\sum_{i\geq j} d_{ij}$ is a differential on $\bigoplus_{i\in I} \idemp_i$, so that the resulting complex is homotopy equivalent to $\one$.
\item $\idempotent_i\otimes \idempotent_j\simeq 0 \cong \idempotent_j\otimes \idempotent_i$ for all $i\neq j$.
\end{enumerate}
Then this gives rise to an idempotent decomposition of $\one$ as follows.  For each convex subset $J\subset I$, let $\idemp_J$ denote the chain complex $\bigoplus_{i\in J}\idemp_i$ with its differential $\sum_{i\leq j\in J}d_{ji}$.  If $J$ is an ideal, then we will also write $\unital_J:=\idemp_J$.  Note that $\unital_I\simeq \one$ by Property (3).  It is a straightforward exercise to verify that $\{\unital_J\}$ is an idempotent decomposition of $\one$ with subquotients given by the $\idemp_J$.
\end{example}

\begin{example}\label{ex:TLex}
Let $\TC\LC_n$ denote Bar-Natan's category $\Mat(\Cob_n)_{\ell}$.  This category is graded additive, monoidal, and has split Grothendieck group isomorphic to an integral form of the Temperley-Lieb algebra.  In \cite{CK12a}, Cooper-Krushkal construct a unital idempotent $P_n\in \KC^-(\TC\LC_n)$ which categorifies the Jones-Wenzl projector (see also \cite{H14a}).  In \cite{CH12}, the author and Ben Cooper extended this to a family of idempotents $P_T\in \KC^-(\TC\LC_n)$, indexed by two-row standard Young tableaux on $n$-boxes.  The categorified Jones-Wenzl projector corresponds then to the one-row tableau.

Let $I=\{T\}$ denote the set of two-row standard tableaux on $n$ boxes.  There is a natural partial order on $I$, given by dominance order; let us in fact order $I$ by the \emph{opposite} of the dominance order, so that the one-row tableau is the unique minimum element.  The main result of \cite{CH12} in fact gives an idempotent decomposition of identity $\{\unital_J\}_{J\subset I}$ with subquotients given by the $P_T$.  
\end{example}

\begin{example}\label{ex:square}
Let $I$ denote the square.  That is, $I=\{a,b,c,d\}$ with partial order $a\leq b,c,\leq d$, with $b$ and $c$ incomparable.  Then the ideals are $\emptyset, \{a\}$, $\{a,b\}$, $\{a,c\}$, $\{a,b,c\}$, and $\{a,b,c,d\}$.  An idempotent decomposition of identity indexed by $I$ is determined by a choice of unital idempotents $\unital_1=\unital_{\{a,b\}}$ and $\unital_2 =\unital_{\{a,c\}}$.  We then have
\[
\unital_\emptyset = 0, \ \ \ \ \ \ \ \ \ \unital_{\{a\}}\cong \unital_1\otimes \unital_2\ \ \ \ \ \ \ \ \ \ \ \ \  \unital_{\{a,b,c\}}\cong \unital_1\vee \unital_2 \ \ \ \ \ \ \ \ \ \ \unital_{\{a,b,c,d\}}\cong \one.
\]
The subquotient idempotents are $\idemp_a\cong \unital_1\otimes \unital_2$, $\idemp_b\cong \unital_1\otimes \unital_2^c$, $\idemp_c\cong \unital_1^c\otimes \unital_2$, and $\idemp_d \cong \unital_1^c\otimes \unital_2^c$.  The decomposition of identity can be depicted diagrammatically by 
\[
\one \simeq \Tot\left(
\begin{diagram}
\unital_1\otimes \unital_2 & \rTo^{[1]} & \unital_1\otimes \unital_2^c \\
\dTo^{[1]} && \dTo^{[1]}\\
\unital_1^c\otimes \unital_2 & \rTo^{[1]} & \unital_1^c\otimes \unital_2^c \\
\end{diagram}
\right).
\]
In order to make this notation precise, it would be necessary to define the notion of a generalized Postnikov system, which is essentially a directed system of objects indexed by a certain partially ordered set.  We will not do this, and instead treat the above notation as a schematic, whose precise meaning can be interpreted in terms of Theorem \ref{thm:generalizedDecomp} below.
\end{example}

Our main theorem in this subsection says that a decomposition of identity can be reconstructed uniquely from its subquotients.  In the following theorem, if $I$ is a partially ordered set, then a \emph{total order} on $I$ is an expression $I=\{i_1,\ldots,i_n\}$ such that $i_k\leq i_\ell$ implies $k\leq \ell$.

\begin{theorem}\label{thm:generalizedDecomp}
Let $I$ be a partially ordered set and $\idemp_i\in \AS$ objects ($i\in I$) such that
\begin{enumerate}
\item $\idemp_i\otimes\idemp_j\cong 0 $ for $i\neq j$.
\item for every total ordering $I= \{i_1,\ldots,i_n\}$, there exists a convolution
\begin{equation}\label{eq:orderedConv}
\one \cong \Tot(\idemp_{i_1}\buildrel [1]\over\rightarrow \cdots \buildrel[1]\over\rightarrow \idemp_{i_n}).
\end{equation}
\end{enumerate}
Then there exists an $I$-indexed decomposition of identity $\{\unital_J,\eta_J\}$ with atomic subquotients the $\idemp_i$.  The idempotents $\unital_J$ determine and are determined by the $\idemp_i$ up to canonical isomorphism.  Furthermore, the idempotents $\idemp_i$ satisfy the semi-orthogonality property
\[
\Homg_\AS(\idemp_i,\idemp_j)=0
\]
unless $i\leq j$.
\end{theorem}
\begin{proof}
We first remark on the uniqueness statement.  Suppose that we have constructed an idempotent decomposition of identity $\{\unital_J,\eta_J\}$ indexed by poset ideals $J\subset I$.  Let $K\subset I$ be a convex set.  Let $J_2$ be the smallest ideal containing $K$, and set $J_1:=J_2\setminus K$.  Then $J_1\subset J_2\subset I$ are ideals, and $K = J_2\setminus J_1$.  We define $\idemp_K:=\DB(\unital_{J_2},\unital_{J_1})$.  We want to ask how does $\idemp_K$ depend on the choice of $J_1\subset J_2$?  So let $J_1'\subset J_2'$ be some other ideals such that $K=J_2'\setminus J_1'$.  The excision property states that $\DB(\unital_{J_2},\unital_{J_1})\cong \DB(\unital_{J_2}',\unital_{J_1}')$, but we must show that this isomorphism is canonical.

The theory of unital idempotents implies that the $\{\unital_J,\eta_J\}$ are unique up to canonical isomorphism.  Proposition \ref{prop:tripleOfRelIdempts} implies that the distinguished triangle which defines $\DB(\unital_{J_2},\unital_{J_1})$ depends uniquely on the idempotents $(\unital_{J_i},\eta_{J_i})$ (and similarly for $J_i'$), up to unique isomorphism of triangles.   From the way that $J_2$ was chosen, we have $J_2\subset J_2'$ and $J_1\subset J_1'$.  Thus, $\unital_{J_i}\leq \unital_{J_i}'$.  Then Proposition \ref{prop:relativeUniqueness} says that the canonical maps $\unital_{J_i'}\rightarrow \unital_{J_i}$ extend to a unique morphism relating the triangles
\[
\DB(\unital_{J_2}',\unital_{J_1}')\rightarrow \unital_{J_2'}\rightarrow \unital_{J_1'}\buildrel[1]\over \rightarrow  \ \ \ \ \ \ \ \ \ \text{and } \ \ \ \ \ \ \ \ \ \ \ \DB(\unital_{J_2},\unital_{J_1})\rightarrow \unital_{J_2}\rightarrow \unital_{J_1}\buildrel[1]\over \rightarrow .
\]
Furthermore, the component from $\DB(\unital_{J_2}',\unital_{J_1}')$ to $\DB(\unital_{J_2},\unital_{J_1})$ is an isomorphism.  This isomorphism, at the end of the day, depends only on the choices of unit maps for $\unital_{J_i}$ and $\unital_{J_i'}$ (and different choices yields the same result up to canonical isomorphism by the same argument).  This takes care of the uniqueness statement.

Now we prove existence.  Let us first define $\idemp_J$.  Let $J\subset I$ be a convex subset, and fix a total ordering $I=\{i_1,\ldots,i_n\}$ such that $J=\{i_p,\ldots,i_q\}$ for $1\leq p\leq q\leq n$.  By hypothesis, we have a convolution of the form (\ref{eq:orderedConv}).   Then Proposition \ref{prop:truncation} gives us a truncation, which we denote by 
\[
\idemp_J = \Tot(\idemp_{i_p}\buildrel [1]\over\rightarrow \cdots \buildrel[1]\over\rightarrow \idemp_{i_q})
\]
Note that $\idemp_J$ satisfies:
\begin{enumerate}
\item[(a)] $\idemp_J$ is in the triangulated hull $\langle \idemp_j\rangle_{j\in J}$.
\item[(b)] $\idemp_J\otimes \idemp_i\cong 0 \cong \idemp_i\otimes \idemp_J$ if $i\not\in J$.
\item[(c)] By reassociation of the convolution (\ref{eq:orderedConv}), $\idemp_J$ fits into a convolution
\[
\one \cong \Tot(\idemp_{i_1}\buildrel [1]\over\rightarrow \cdots \buildrel [1]\over\rightarrow \idemp_{i_{p-1}} \buildrel [1]\over\rightarrow  \idemp_J \buildrel [1]\over\rightarrow \idemp_{i_{q+1}} \buildrel [1]\over\rightarrow  \cdots \buildrel [1]\over\rightarrow  \idemp_{i_n}).
\]
If $C\in \TS$ is any object in the triangulated hull of the $\idemp_j$ ($j\in J$), then $C\otimes(-)$ and $(-)\otimes C$ annihilate all terms of the above except for $\idemp_J$, hence tensoring this convolution with $C$ yields
\[
\idemp_J\otimes C\cong C \cong C\otimes \idemp_J
\]
for all $C\in \langle \idemp_j\rangle_{j\in J}$.
\end{enumerate}
The properties (a),(b),(c) uniquely characterize $\idemp_J$ up to isomorphism.  Properties (a) and (c) imply that $\idemp_J^{\otimes 2}\cong \idemp_J$.  We now make the following additional observations:
\begin{enumerate}
\item[(d)] If $J,K\subset I$ are disjoint convex sets, then $\idemp_J\otimes \idemp_K\cong 0$.
\item[(e)] If $J\subset I$ is an ideal, then reassociation of (\ref{eq:orderedConv}) gives
\[
\one\cong \Tot(\idemp_J\buildrel [1]\over \rightarrow \idemp_{I\setminus J}),
\]
which simply means there is a distinguished triangle
\[
\idemp_{I\setminus J}\rightarrow \one \rightarrow \idemp_J \rightarrow \idemp_{I\setminus J}[1].
\]
But $\idemp_J\ \otimes \idemp_{I\setminus J}\cong 0 \cong \idemp_{I\setminus J} \otimes \idemp_J$ by (d), hence $\unital_J:=\idemp_J$ has the structure of a unital idempotent with complementary idempotent $\unital_J^c = \idemp_{I\setminus J}$.
\item[(f)] If $K_1\subset K_2\subset I$ are ideals, then $\unital_{K_1}$ is orthogonal to $\unital_{K_2}^c$ by (d), which implies that $\unital_{K_1}\leq \unital_{K_2}$.
\item[(g)]  We have $\idemp_J\cong \DB(\unital_{K_2},\unital_{K_1})$, by Proposition \ref{prop:subquotientCharacterization}.
\end{enumerate}

Now, let us show that properties (1) and (2) of Definition \ref{def:decompOfOne} are satisfied.  Let $J,K\subset I$ be ideals.  From property (g), there exists a distinguished triangle
\[
\idemp_{J\setminus K} \rightarrow \unital_J\rightarrow \unital_{J\cap K}\rightarrow \idemp_{J\setminus K} [1],
\]
Tensoring with $\unital_K$ yields a distinguished triangle
\[
0\rightarrow \unital_J\otimes \unital_K \rightarrow \unital_{J\cap K}\otimes \unital_K\rightarrow 0.
\]
The first term is zero by property (d).  The third term is isomorphic to $\unital_{J\cap K}$ since $\unital_{J\cap K}\leq \unital_K$ by property (f).  This shows that $\unital_J\otimes \unital_K\cong \unital_{J\cap K}\cong \unital_J\otimes \unital_K$.  In particular, the idempotents $\unital_J$ commute with one another.

We wish to show that $\unital_{J\cup K}\cong \unital_J\vee \unital_K$.   First, observe that $\unital_J\leq \unital_{J\cup K}$ and $\unital_K\leq \unital_{J\cup K}$, hence $\unital_J\vee \unital_K\leq \unital_{J\cup K}$.  On the other hand, the complementary idempotent satisfies
\[
(\unital_J\vee \unital_K)^c\cong \unital_J^c\otimes \unital_K^c \cong \idemp_{I\setminus J}\otimes \idemp_{I\setminus K}.
\]
This last expression clearly is annihilated by $\idemp_i$ whenever $i\in J\cup K$ (the first factor kills the $\idemp_i$ with $i\in J$, and the second factor kills the $\idemp_i$ with $i\in K$).  Thus, $(\unital_J\vee \unital_K)^c$ is annihilated by $\unital_{J\cup K}$.  It follows that the complementary idempotent $\unital_J\vee \unital_K$ fixes $\unital_{J\cup K}$.  That is, $\unital_{J\cup K}\leq \unital_J\vee \unital_K$.  This proves that $\unital_{J\cup K}\cong \unital_1\vee\unital_2$.

Finally, the semi-orthogonality property follows along the same lines as in Proposition \ref{prop:subquotientProps}.  If $i\not< j$, then we may find a ideals $J$ containing $j$ but not containing $i$.  It follows that $\unital_J\otimes \idemp_j\cong \idemp_j$ and $\unital_J^c\otimes \idemp_i \cong \idemp_i$.  Then the relation $\Homg_{\AS}(\idemp_i,\idemp_j)\cong 0$ follows from Proposition \ref{thm-semiOrtho}.
\end{proof}

Finally, we have the following:

\begin{proposition}\label{prop:tensorProductResolution}
Let $\{\unital_J\}_I$ and $\{\unital_J'\}_{I'}$ be idempotent decompositions of $\one$ indexed by partially ordered sets $I$ and $I'$, with atomic subquotients $\idemp_i$, $\idemp_{i'}'$, respectively.  If the $\idemp_i$ and $\idemp_{i'}'$ commute for all $i\in I$, $i',\in I'$, then $\{\unital_J\otimes \unital_{J'}'\}_{I\times I'}$ extends to an idempotent decomposition of identity with atomic subquotients $\idemp_i\otimes \idemp_{i'}'$.  Here the indexing set $I\times J$ has the product partial order: $(i,j)\leq (i',j')$ if $i\leq i'$ and $j\leq j'$. 
\end{proposition}
\begin{proof}
Let regard $I,I'$, and $I\times I'$ as topological spaces with open sets given by the ideals.  An easy exercise shows that $I\times I'$ has the product topology.  By a straightforward induction, the fact that the subquotients commute implies (using Proposition \ref{prop:commuteWithidemp}), that the idempotents $\unital_J$ and $\unital_{J'}'$ commute for all open sets $J\subset I$, $J'\subset I'$.   Thus, Proposition \ref{prop:suprema} implies that $\unital_J\otimes \unital_{J'}'\cong \unital_J\wedge \unital_{J'}'$ is a unital idempotent.  

The subsets $J\times J'\subset I\times I'$ form a basis for the topology on $I\times I'$, hence we may extend to a collection of idempotents $\unital''_K$, indexed by open subsets $K\subset I\times I'$.  We leave it to the reader to verify that the idempotents so obtained satisfy the requirements of Definition \ref{def:decompOfOne}.
\end{proof}

\section{Idempotents and Tate cohomology}
\label{sec:tate}

Let $\k$ be a field and $H$ a finite dimensional Hopf algebra over $\k$.  Below, we will abbreviate $\Homg_{\KC(H\modules)}(-,-)$ simply by $\Homg_H(-,-)$.   Let $\e:\PB\rightarrow \k$ denote a resolution of the trivial $H$-module $\k$ by projective $H$-modules.  Then $(\PB,\e)$ is a counital idempotent in $\KC^-(H\modules)$.  By definition, the \emph{cohomology of $H$} is $\Ext_H(\k,\k)$, which can be computed in terms of $\PB$ as
\[
\Ext_H(\k,\k) \cong \Endg_H(\PB) \cong \Homg_H(\PB,\k).
\]

The \emph{homology of $H$} is by definition $\Tor_H(\k,\k)$, which can be computed as the homology of $\PB\otimes_\k \k$ or, equivalently
\[
\Tor_H(\k,\k) \cong \Homg_H(\k,\PB).
\]

The Tate cohomology is obtained by patching together the cohomology and homology of $H$.  To be precise, recall that if $M$ is an $H$-module and $S:H\rightarrow H$ denotes the antipode, then the linear dual $M^\star:=\Hom_\k(M,\k)$ is an $H$-module via $x\phi:m\mapsto \phi(S(x)m)$ for all $x\in H$, $\phi\in M^\ast$, $m\in M$.  The dual map $\e^\star:\k\rightarrow \PB^\star$ is an injective resolution of $\k$.  Finite dimensional Hopf algebras are Frobenius, which implies that injectives coincide with projectives, hence $\PB^\star$ could also be thought of as a projective coresolution of $\k$.  Then we define an object $\Tate:=\Cone(\PB\buildrel \e^\star\circ \e\over \longrightarrow \PB^\star)$. The \emph{Tate cohomology of $H$} is by definition the homology of $\Tate$.  Sometimes this is denoted by $\widehat{\Ext}_H(\k,\k)$.  In other words:
\[
\widehat{\Ext}_H(\k,\k) = \Homg_H(\k,\Tate).
\] 
Applying the functor $\Homg_H(\k,-)$ to the distinguished triangle $\PB\rightarrow \PB^\star\rightarrow \Tate\rightarrow \PB[1]$ gives rise to a long exact sequence
\[
\Homg_H(\k,\PB)\rightarrow \Homg_H(\k,\PB^\star)\rightarrow \Homg_H(\k,\Tate)\rightarrow \Homg_H(\k,\PB)[1].
\]
The first term is the homology of $H$, the second is the cohomology of $H$, and the third is the Tate cohomology of $H$.

One of the main goals of this section is to prove that the complex $\Tate$ has the structure of a unital idempotent in $\KC(H\modules)$, hence the Tate cohomology can be computed as 
\[
\widehat{\Ext}_H(\k,\k) \cong \Endg_H(\Tate).
\]
This implies that the Tate cohomology has the structure of a graded commutative ring, which is not entirely obvious from the definition.  
\begin{remark}
One has to be careful when tensoring together unbounded complexes.  There are two ways to do this, one using direct sum and the other using direct product (See equations (\ref{eq:tensorProdSum}) and (\ref{eq:tensorProdProd})).  The former is denoted by $\otimes$, and the latter by $\cotimes$.  As it turns out, $\Tate[-1]$ is a unital idempotent with respect to $\cotimes$.  Thus, the Tate cohomology could also be computed as 
\[
\widehat{\Ext}_H(\k,\k) \cong \Endg_H(\Tate)\cong \Homg_H(\Tate[-1],\k)\cong \Homg_H(\k,\Tate^\star[1])
\]
This is consistent with what we wrote above since $\Tate^\star = \Cone(\PB\rightarrow \PB^\star)^\star \cong \Cone({\PB^\star}^\star\rightarrow \PB^\star)[-1]\cong \Tate[-1]$, because of the shifts involved in forming the mapping cone.
\end{remark}

\begin{remark}
One can define Tate cohomology in the more general setting of strongly Gorenstein rings \cite{BuchweitzPreprint}.  In this more general setting, the duality functor is less well-behaved.  The results of \S \ref{subsec:generalTate} apply to this more general setting, but the results of \ref{subsec:duals} apply only to the special case in which there is a ``nice'' notion of duality.
\end{remark}

\subsection{Generalized Tate cohomology}
\label{subsec:generalTate}
Let $(\counital_i,\unital_i)$ ($i=1,2$) be two pairs of complementary idempotents; they fit into distinguished triangles of the form
\begin{equation}\label{eq:T3triang1}
\counital_1 \buildrel \e_1\over \longrightarrow \one \buildrel \eta_1\over \longrightarrow \unital_1\buildrel \d_1\over\longrightarrow \counital_1[1]
\end{equation}
and
\begin{equation}\label{eq:T3triang2}
 \one \buildrel \eta_2\over \longrightarrow \unital_2\buildrel \d_2\over\longrightarrow \counital_2[1] \buildrel -\e_2[1]\over \longrightarrow \one[1].
\end{equation}
Choose an object $\Tate:= \Cone(\eta_2\circ \e_1)$ and a distinguished triangle
\begin{equation}\label{eq:TateTriang0}
\begin{diagram}[small]
\counital_1 & \rTo^{\eta_2\circ \e_1} & \unital_2 & \rTo^{\nu_2} & \Tate & \rTo^{\theta_1} & \counital_1[1]
\end{diagram}.
\end{equation}
From the axioms of triangulated categories (see axiom T3 in \cite{MayTraces}) there exists a distinguished triangle
\begin{equation}\label{eq:TateTriang1}
\begin{diagram}[small]
\counital_2 & \rTo^{\eta_1\circ \e_2} & \unital_1 & \rTo^{\nu_1} & \Tate & \rTo^{\theta_2} & \counital_2[1]
\end{diagram}
\end{equation}
such that the following diagram commutes:
\begin{equation}\label{eq:TateDiagram1}
\begin{minipage}{5.1in}
\begin{tikzpicture}
\tikzstyle{every node}=[font=\small]
\node (a) at (0,0) {$\counital_1$};
\node (b) at (4,0) {$\unital_2$};
\node (c) at (8,0) {$\counital_2[1]$};
\node (d) at (12,0) {$\unital_1[1]$};
\node (e) at (2,-2) {$\one$};
\node (f) at (6,-2) {$\Tate$};
\node (g) at (10,-2) {$\one[1]$};
\node (h) at (4,-4) {$\unital_1$};
\node (i) at (8,-4) {$\counital_1[1]$};\\
\path[->,>=stealth',shorten >=1pt,auto,node distance=1.8cm, thick]
(a) edge node {$\e_1$} (e)
(e) edge node {$\eta_1$} (h)
(h) edge [bend right, looseness =1,out=-45,in=225] node[below] {$\d_1$} (i)
(i) edge node[label={[shift={(.7,-.7)}]$\e_1[1]$}] {}  (g)
(g) edge node[label={[shift={(.7,-.7)}]$\eta_1[1]$}] {}  (d)
(e) edge node[below right] {$\eta_2$} (b)
(b) edge [bend left, looseness=1,out=45,in=135]  node {$\d_2$} (c)
(c) edge node {$-\e_2[1]$} (g)
(c) edge [bend left, looseness=1,out=45,in=135] node {$(\eta_1\circ \e_2)[1]$} (d)
(b) edge node {$\nu_2$} (f)
(f) edge node {$\theta_1$} (i)
(a) edge [bend left, looseness=1,out=45,in=135]  node {$\eta_2\circ \e_1$} (b)
;
\path[dashed,->,>=stealth',shorten >=1pt,auto,node distance=1.8cm, thick]
(h) edge node[below right] {$\nu_1$} (f)
(f) edge node[below right] {$\theta_2$} (c);
\end{tikzpicture}
\end{minipage}
\end{equation}

In particular $\Tate\simeq \Cone(\eta_1\circ \e_2)\simeq \Cone(\eta_2\circ \e_1)$.

\begin{definition}\label{def:TateObject}
Let $(\counital_i,\unital_i)$ ($i=1,2$) be pairs of complementary idempotents.  Define $\Tate$ to be the object $\Tate\simeq \Cone(\eta_1\circ \e_2)\simeq \Cone(\eta_2\circ \e_1)$ described above.  We refer to $\Tate$ as the Tate object associated to $(\counital_2,\unital_1)$ or, equivalently, $(\counital_1,\unital_2)$, and we will write $\Tate\simeq \Tate(\counital_2,\unital_1)\simeq \Tate(\counital_1,\unital_2)$.

 The \emph{Tate cohomology} of the pair $(\counital_1,\unital_2)$ or, equivalently $(\counital_2,\unital_1)$, is  the graded algebra $\Homg_\AS(\Tate,\Tate)$.  We will usually write this as $H^\bullet(\counital_1,\unital_2)$ or, equivalently, $H^\bullet(\counital_2,\unital_1)$.
\end{definition}

\begin{example}\label{ex:Tate}
We have the following (somewhat degenerate) examples of Tate objects.  Let $(\counital,\unital)$ be a pair of complementary idempotents.   Then
\begin{enumerate}
\item $\Tate(\counital, \one)\simeq \unital$.
\item $\Tate(\one,\unital)\simeq \counital[1]$.
\item $\Tate(\counital,\unital)\simeq \unital\oplus \counital$.
\end{enumerate}
\end{example}

Note that $\Tate$ is only unique up to (possibly non-canonical) isomorphism.  However, under certain assumptions, $\Tate$ is a categorical idempotent, hence satisfies a stronger form of uniqueness by Theorem \ref{thm-uniquenessOfIdempts}.


\begin{theorem}\label{thm:fundThmOfTateIdempt}
Retain notation as above.  In addition to the distinguished triangles (\ref{eq:TateTriang0}) and (\ref{eq:TateTriang1}), the object $\Tate$ fits into distinguished triangles:
\begin{equation}\label{eq:TateTriang2}
\begin{diagram}[small]
(\counital_1\oplus \counital_2) & \rTo^{\matrix{\e_1 & \e_2}} &  \one & \rTo & \Tate &\rTo & (\counital_1\oplus \counital_2)[1]
\end{diagram}
\end{equation}
and
\begin{equation}\label{eq:TateTriang3}
\begin{diagram}[small]
\Tate[-1] & \rTo &\one & \rTo^{\matrix{\eta_1\\ \eta_2}} & (\unital_1\oplus \unital_2) & \rTo & \Tate.
\end{diagram}
\end{equation}
Furthermore, if $\counital_1\otimes \counital_2\simeq 0 \simeq \counital_2\otimes \counital_1$, then $\Tate\simeq \unital_1\otimes \unital_2$ has the structure of a unital idempotent with complement $\counital_1\oplus \counital_2$.

On the other hand if $\unital_1\otimes \unital_2\simeq 0 \simeq \unital_2\otimes \unital_1$,  then $\Tate[-1]\simeq \counital_1\otimes \counital_2$ has the structure of a counital idempotent with complement $\unital_1\oplus \unital_2$.
\end{theorem}

\begin{proof}
This is an easy consequence of the $3\times 3$-lemma (see  Lemma 2.6 in \cite{MayTraces}).  First, consider the following diagram
\[
\begin{tikzpicture}
\tikzstyle{every node}=[font=\small]
\node (aa) at (0,0) {$\counital_1$};
\node (ba) at (2.5,0) {$\counital_1$};
\node (ca) at (5,0) {$0$};
\node (da) at (7.5,0) {$\counital_1[1]$};
\node (ab) at (0,-2.5) {$\counital_1\oplus \counital_2$};
\node (bb) at (2.5,-2.5) {$\one$};
\node (cb) at (5,-2.5) {$X$};
\node (db) at (7.5,-2.5) {$\counital_1[1]\oplus \counital_2[1]$};
\node (ac) at (0,-5) {$\counital_2$};
\node (bc) at (2.5,-5) {$\unital_1$};
\node (cc) at (5,-5) {$Y$};
\node (dc) at (7.5,-5) {$\counital_2[1]$};
\node (ad) at (0,-7.5) {$\counital_1[1]$};
\node (bd) at (2.5,-7.5) {$\counital_1[1]$};
\node (cd) at (5,-7.5) {$0$};
\node (dd) at (7.5,-7.5) {$\counital_1[2]$}
;\\
\path[->,>=stealth',shorten >=1pt,auto,node distance=1.8cm, thick]
(aa) edge node {$\Id$} (ba)
(ba) edge node {} (ca)
(ca) edge node {} (da)
(ab) edge node {$\matrix{\e_1 & \e_2}$} (bb)
(bb) edge node {} (cb)
(cb) edge node {} (db)
(ad) edge node {$\Id$} (bd)
(bd) edge node {} (cd)
(cd) edge node {} (dd)
(aa) edge node {$\matrix{\Id\\ 0}$} (ab)
(ab) edge node {$\matrix{0 & \Id}$} (ac)
(ac) edge node {$0$} (ad)
(ba) edge node {$\e_1$} (bb)
(bb) edge node {$\eta_1$} (bc)
(bc) edge node {$\d_1$} (bd)
(da) edge node {$\matrix{\Id\\ 0}$} (db)
(db) edge node {$\matrix{0 & \Id}$} (dc)
(dc) edge node {$0$} (dd)
;
\path[dashed,->,>=stealth',shorten >=1pt,auto,node distance=1.8cm, thick]
(ac) edge node {$f$} (bc)
(bc) edge node {} (cc)
(cc) edge node {} (dc)
(ca) edge node {} (cb)
(cb) edge node {$g$} (cc)
(cc) edge node {} (cd)
;
\end{tikzpicture}
\]
The first two rows and the first two columns are distinguished, and the square on the top left commutes.  The $3\times 3$ lemma guarantees that there exist morphisms corresponding to the dashed arrows above, so that the third row and column are distinguished, and that every square commutes except for the square on the bottom right, which commutes up to sign.

Commutativity of the middle square on the left forces $f=\eta_1\circ \e_2$.  Thus $Y\simeq \Tate$, by uniqueness of mapping cones. The third column being distinguished implies that $g:X\rightarrow Y\simeq \Tate$ is an isomorphism; the second row is then the distinguished triangle (\ref{eq:TateTriang2}) we were seeking.  A similar argument constructs the distinguished triangle (\ref{eq:TateTriang3}).

Now, assume that $\counital_1\otimes \counital_2\simeq 0 \simeq \counital_2\otimes \counital_1$.  Tensoring the distinguished triangle (\ref{eq:TateTriang1}) on the left with $\counital_1$ gives a distinguished triangle $0\rightarrow 0 \rightarrow \counital_1\otimes \Tate\rightarrow 0$, which implies that $\counital_1\otimes \Tate\simeq 0$.  Similar arguments show that
\[
0\simeq \counital_2\otimes \Tate\simeq \Tate\otimes \counital_1\simeq \Tate\otimes \counital_2
\]
This, together with the distinguished triangle (\ref{eq:TateTriang2}), implies that $(\Tate, \counital_1\oplus\counital_2)$ is a pair of complementary idempotents.  To show that $\unital_1\otimes \unital_2\simeq \Tate$, first observe that $\Tate$ absorbs the unital idempotents $\unital_1,\unital_2$, since it annihilates their complementary idempotents $\counital_1,\counital_2$.  Now, tensor the triangle (\ref{eq:TateTriang2}) on the left by $\unital_1$ and on the right by $\unital_2$, obtaining a distinguished triangle
\[
0\rightarrow \unital_1\otimes \unital_2 \rightarrow \Tate \rightarrow 0.
\]
This implies that $\Tate\simeq \unital_1\otimes \unital_2$, as claimed.  An entirely similar argument proves the final statement, starting from the distinguished triangle (\ref{eq:TateTriang3}).
\end{proof}

\subsection{Duals}
\label{subsec:duals}
Let $\AS$ be a triangulated monoidal category.  In this section we consider the relationship between categorical idempotents in $\AS$ and their duals (assuming some sort of duality is present).  Perhaps the most general setting in which these results can be stated is that of rigid triangulated monoidal categories \cite{MayTraces}.  Since categorical idempotents tend to be ``large'' (e.g.~infinite complexes), the notion of dual is not well-behaved.  In order to avoid the subtleties involved, from this point on we restrict to examples in which our triangulated categories are full subcategories of the homotopy category $\KC(\CS)$, for an appropriate additive category $\CS$.

 In typical examples $\CS$ is an additive monoidal category, and we have a pair of complementary idempotents $(\counital,\unital)$ in $\KC^-(\CS)$.  Their duals (if defined) will be a pair $(\counital^\star,\unital^\star)$ of complementary idempotents $\KC^+(\CS)$.  In order to compare the idempotents and their duals, we must pass to the category $\KC(\CS)$ of unbounded complexes.  There are two different ways to take the tensor product of unbounded complexes, which we now explain.

Let $\CS^\oplus$ and $\CS^\Pi$ denote the categories obtained from $\CS$ by formally adjoining countable sums and products to $\CS$, respectively.  An object of $\CS^\oplus$ is a formal sequence $(A_0,A_1,\ldots,)$, where $A_i$ are objects of $\CS$.  A morphism in $\CS$ is an $\N\times \N$-matrix of morphisms in $\CS$ which is \emph{column-finite} (in each column, at most finitely many entries are nonzero).  To define $\CS^\Pi$ replace the column-finiteness condition with row-finiteness.  Inside $\CS^\oplus$ (respectively $\CS^\Pi$) we have $(A_0,A_1,\ldots) \cong \bigoplus_{i\in \N} A_i$  (respectively $(A_0,A_1,\ldots) \cong \prod_{i\in \N} A_i$).  Note that $\CS$ includes as a full subcategory of $\CS^\oplus$ and $\CS^\Pi$, via $A\mapsto (A,0,0,\ldots )$.

There is a monoidal structure on $\KC(\CS^\oplus)$, with $\otimes$ being defined in the usual way (Equation (\ref{eq:tensorProdSum})).  Replacing the direct sum with a direct product defines the other tensor product $A\cotimes B$ (Equation (\ref{eq:tensorProdProd})).  Then $\KC(\CS^\Pi)$ has the structure of a monoidal category with monoid $\cotimes$.

Below, let $\CS$ be an additive, rigid monoidal category with duality $(-)^\star:\CS\rightarrow \CS^{\op}$.  In particular, there are isomorphisms
\[
\Hom_{\CS}(a,b)\cong \Hom_{\CS}(a\otimes b^\star,\one)\cong \Hom_{\CS}(\one,b\otimes a^\star)
\]
and
\[
(a\otimes b)^\star \cong b^\star\otimes a^\star
\]
which are natural in $a,b\in \CS$.   The duality functor extends to complexes.  By convention:
\[
(A^\star)_{i} = (A_{-i})^\star \ \ \ \ \text{with differential} \ \ \ \ d_{(A^\star)_i} := (-1)^{i+1}(d_{A_{-i-1}}^\star)
\]
It is an easy exercise to show that for any complexes $A,B\in \KC(\CS)$, we have
\begin{equation}\label{eq:dualityForComplexes}
\Hom_{\KC(\CS)}(A,B) \cong \Hom_{\KC(\CS^\oplus)} (A\otimes B^\star, \one)\cong \Hom_{\KC(\CS^\Pi)}(\one,B\cotimes A^\star)
\end{equation}

As another easy exercise, if $f:A\rightarrow B$ is a chain map, then
\begin{equation}\label{eq:dualOfCone}
\Cone(f)^\star\cong \Cone(f^\star)[1].
\end{equation}

In fact, $(-)^\star$ gives a pair of inverse functors $\KC^\pm(\CS)\rightarrow \KC^\mp(\CS)$ and $\KC(\CS^\oplus)\leftrightarrow \KC(\CS^\Pi)$.    The most important properties of dual idempotents are collected below:

\begin{proposition}\label{prop:dualIdempts}
Let $(\counital,\unital)$ be a pair of complementary idempotents in $\KC^-(\CS)$.  Then $(\counital^\star,\unital^\star)$ is a pair of complementary idempotents in $\KC^+(\CS)$.  Moreover,
\begin{equation}\label{eq:counitalAndDual}
\counital\otimes \unital^\star \simeq 0 \simeq \unital^\star \otimes \counital \ \ \ \ \ \ \ \in \ \ \  \KC(\CS^\oplus)
\end{equation}
\begin{equation}\label{eq:counitalAndDual2}
\counital^\star\cotimes \unital \simeq 0 \simeq \unital \cotimes \counital^\star \ \ \ \ \ \  \in \ \ \  \KC(\CS^\Pi)
\end{equation}
\end{proposition}
\begin{proof}
We have a distinguished triangle of the form
\[
\counital\buildrel\e\over \longrightarrow \one \buildrel\eta\over \longrightarrow \unital \rightarrow \counital[1]
\]
in $\KC^-(\CS)$.  Taking duals gives a distinguished triangle
\[
\unital^\star\buildrel\eta^\star\over \longrightarrow \one \buildrel \e^\star\over \longrightarrow \counital^\star\rightarrow \unital^\star[1].
\]
Furthermore, $\unital^\star\otimes \counital^\star \cong (\counital\otimes \unital)^\star\simeq 0$, and $\counital^\star\otimes \unital^\star\cong (\unital\otimes \counital)^\star\simeq 0$.  Thus, $(\counital^\star,\unital^\star)$ is a pair of complementary idempotents in $\KC^+(\CS)$.

Now, $\e\otimes \Id_{\counital} : \counital\otimes \counital\rightarrow \one\otimes \counital$ is a homotopy equivalence by projector absorbing.  Similarly, $\Id_{\unital^\star}\otimes \eta^\star: \unital^\star\otimes \unital^\star\rightarrow \unital^\star\otimes \one$ is a homotopy equivalence.  The tensor product of these maps is a homotopy equivalence
\[
\Id_{\unital^\star}\otimes \eta^\star\otimes \e\otimes \Id_\counital : \unital^\star\otimes \unital^\star\otimes \counital\otimes \counital \buildrel \simeq \over \rightarrow \unital^\star\otimes \one\otimes \one \otimes \counital.
\]
On the other hand this map is homotopic to zero, since
\[
[\eta^\star\otimes \e] \in \Homg(\unital^\star\otimes \counital,\one) \cong \Homg(\counital,\unital)\cong 0
\]
by semi-orthogonality of idempotents.  This proves $\unital^\star\otimes \counital\simeq 0$.  The remaining equivalences are proven similarly.
\end{proof}

\begin{theorem}\label{theorem:dualIdempotents}
Let $(\counital,\unital)$ be a pair of complementary idempotents in $\KC^-(\CS)$, and let $(\counital^\star,\unital^\star)$ be the dual pair in $\KC^+(\CS)$.  Let $\Tate\in \KC(\CS)$ be the Tate object associated to $(\unital^\star,\unital)$ (Definition \ref{def:TateObject}).  Then
\begin{enumerate}
\item $\Tate$ has the structure of a unital idempotent in $\KC(\CS^\oplus)$ with complementary idempotent $\counital\oplus \unital^\star$.
\item $\Tate[-1]$ has the structure of a counital idempotent in $\KC(\CS^\Pi)$ with complementary idempotent $\unital\oplus \counital^\star$. 
\end{enumerate}
  Furthermore, we have
\begin{eqnarray}
\Homg(\Tate,\Tate)
&\simeq & \Homg(\one,\Tate) \label{eq:Tate1}\\
&\simeq & \Homg(\Tate[-1],\one)\label{eq:Tate2}\\
&\simeq &  \Homg(\unital,\counital[1])\label{eq:Tate3}\\
&\simeq & \Homg(\counital^\star,\unital^\star[1])\label{eq:Tate4}
\end{eqnarray}
Each of these isomorphisms commutes with the action of $\Homg(\unital,\unital)\otimes\Homg(\counital,\counital)$.  The isomorphism $\Homg(\Tate,\Tate)\rightarrow \Homg(\unital,\counital[1])$ sends $[\Id_\Tate]\mapsto [\d]$, where $\d:\unital\rightarrow \counital[1]$ is the connecting morphism.  Thus, $\Homg(\unital,\counital[1])$ has the structure of a unital algebra with unit $[\d]$.
\end{theorem}
\begin{proof}
Statements (1) and (2) are immediate consequences of Theorem \ref{thm:fundThmOfTateIdempt} together with Proposition \ref{prop:dualIdempts}.  The equivalences (\ref{eq:Tate1}) and (\ref{eq:Tate2})  then follow from general arguments (see Corollary \ref{cor-endRings}).  We now prove (\ref{eq:Tate3}).  

Note that $\Homg(\unital,\unital)\cong \Homg(\unital^\star,\unital^\star)$ (the isomorphism given by the functor $(-)^\star$).  We may think of $\Tate$ as a counital idempotent relative to $\unital^\star\oplus \counital$, hence we have an algebra map $\Homg(\unital^\star,\unital^\star)\otimes \Homg(\counital,\counital)\rightarrow \Homg(\Tate,\Tate)$.  All of the isomorphisms below will commute with this action by statement (6) of Proposition \ref{prop:superActions}.

 Recall the diagram (\ref{eq:TateDiagram1}), with $\unital_1=\unital$ and $\unital_2=\counital^\star$. The map $\nu:\unital\rightarrow \Tate$ fits into a distinguished triangle $\unital^\star\rightarrow \unital \rightarrow \Tate\rightarrow \unital^\star$, which realizes $\Tate$ as a unital idempotent relative to $\unital$ (since $\Tate$ annihilates $\unital^\star$ in $\KC(\CS^\oplus)$).  Thus, precomposing with $\nu$ gives an isomorphism $\Homg(\Tate,\Tate)\rightarrow \Homg(\unital,\Tate)$ sending $[\Id_\Tate]\mapsto [\nu]$.  Now, we may also think of $\Tate[-1]$ as a counital idempotent relative to $\counital$, with structure map $\theta[-1]:\Tate[-1]\rightarrow \counital$.  We claim that post-composing with $\theta$ gives an isomorphism
 \[
 \Homg(\unital,\Tate)\rightarrow \Homg(\unital,\counital[1])
 \]
 By construction, this isomorphism will send $[\Id]\mapsto [\theta\circ \nu]$, which equals $[\d]$ by commutativity of the top middle triangle in Diagram \ref{eq:TateDiagram1}. 
 
 Observe that $\Tate$ fits into a distinguished triangle $\counital\rightarrow \counital^\star\rightarrow \Tate\rightarrow \counital[1]$ in which the last map is $\theta$.  We have
 \[
 \Homg(\unital,\counital^\star)\cong \Homg(\unital\otimes \counital,\one)\cong 0.
 \]
 Thus, applying $\Homg(\unital,-)$ to the distinguished triangle $\counital\rightarrow \counital^\star\rightarrow \Tate\rightarrow \counital[1]$ gives a long exact sequence
 \[
\Homg(\unital,\counital)\rightarrow \Homg(\unital,\counital^\star)\rightarrow \Homg(\unital,\Tate)\buildrel \theta_\ast\over\rightarrow \Homg(\unital,\counital[1])
 \]
 The second term is zero, hence the third map is an isomorphism, as claimed.  This proves (3).  A similar argument (or applying the funtor $(-)^\star$) gives (4).
 \end{proof}

\printbibliography

\end{document}